\documentclass{article}
\usepackage{graphicx} 
\usepackage[noblocks]{authblk}
\usepackage{amsmath}
\usepackage{amsthm}
\usepackage{amsfonts}
\usepackage{subcaption} 
\usepackage{float}
\usepackage{algorithm}
\usepackage{algpseudocode}
\usepackage{xcolor}
\usepackage{booktabs}
\usepackage{array}
\usepackage{youngtab}
\usepackage{xparse}
\usepackage{longtable}
\usepackage[colorlinks=true, linkcolor=blue, citecolor=blue, urlcolor=blue]{hyperref}
\usepackage{array,longtable,booktabs}
\NewDocumentCommand{\swilf}{ o }{%
  \overset{\mathrm{s}\IfValueT{#1}{(#1)}}{\sim}%
}
\NewDocumentCommand{\wilf}{ o }{%
  \overset{\IfValueT{#1}{(#1)}}{\sim}%
}

\NewDocumentCommand{\PatternFormat}{m m}{{#1}_{#2}}
\NewDocumentCommand{\PatternSetFormat}{m}{\{#1\}}
\newcommand{\PatternSeparator}{,}

\ExplSyntaxOn
\tl_new:N \l__paper_pat_values_tl
\tl_new:N \l__paper_pat_colors_tl
\bool_new:N \l__paper_pat_first_bool

\cs_new_protected:Npn \__paper_pat_parse:w #1/#2\q_stop
  {
    \tl_set:Nn \l__paper_pat_values_tl {#1}
    \tl_set:Nn \l__paper_pat_colors_tl {#2}
    \tl_trim_spaces:N \l__paper_pat_values_tl
    \tl_trim_spaces:N \l__paper_pat_colors_tl
    \exp_args:NVV \PatternFormat
      \l__paper_pat_values_tl \l__paper_pat_colors_tl
  }

\NewDocumentCommand{\pat}{m}
  {
    \ensuremath
      {
        \group_begin:
        \__paper_pat_parse:w #1\q_stop
        \group_end:
      }
  }

\cs_new_protected:Npn \__paper_pat_list:n #1
  {
    \bool_set_true:N \l__paper_pat_first_bool
    \clist_map_inline:nn {#1}
      {
        \bool_if:NTF \l__paper_pat_first_bool
          { \bool_set_false:N \l__paper_pat_first_bool }
          { \PatternSeparator }
        \pat{##1}
      }
  }

\NewDocumentCommand{\patset}{s m}
  {
    \ensuremath
      {
        \group_begin:
        \PatternSetFormat
          { \IfBooleanTF{#1}{#2}{\__paper_pat_list:n {#2}} }
        \group_end:
      }
  }
\ExplSyntaxOff

\newcommand{\NE}{\operatorname{NE}}
\newcommand{\NW}{\operatorname{NW}}
\newcommand{\SE}{\operatorname{SE}}
\newcommand{\SW}{\operatorname{SW}}
\newcommand{\rowSE}{\operatorname{rowSE}}
\newcommand{\rowSW}{\operatorname{rowSW}}

\definecolor{rookR}{HTML}{B43138}
\definecolor{rookB}{HTML}{2461AD}
\definecolor{marked}{HTML}{D9DDE1}

\newcommand{\drawstackplacement}[3]{%
  \foreach \height [count=\col] in {#1} {
    \ifnum\col=#3\relax
      \path[fill=black!10]
        ({\col-1},0) rectangle (\col,\height);
    \fi
    \foreach \row in {1,...,\height} {
      \draw[boardcell]
        ({\col-1},{\row-1}) rectangle ++(1,1);
    }
  }
  \foreach \row [count=\col] in {#2} {
    \node at ({\col-0.5},{\row-0.5}) {$1$};
  }
}
\usepackage[utf8]{inputenc}
\usepackage{xcolor}
\usepackage{array}

\usepackage[normalem]{ulem} 
\usepackage{todonotes}
\usepackage{tikz}
\usetikzlibrary{patterns,arrows.meta}

\newcommand{\stkout}[1]{\ifmmode\text{\sout{\ensuremath{#1}}}\else\sout{#1}\fi}

\newtheorem{lemma}{Lemma}[section]
\newtheorem{theorem}[lemma]{Theorem}
\newtheorem{corollary}[lemma]{Corollary}
\newtheorem{proposition}[lemma]{Proposition}
\newtheorem{observation}[lemma]{Observation}

\theoremstyle{definition}

\begin{document}

\title{Shape-Wilf Equivalence for Colored Patterns}
\author{Tigran Arsenyan}
\author{V\'it Jel\'inek}

\affil{Computer Science Institute, Faculty of Mathematics and Physics, Charles University, Prague}
\date{}
\maketitle

\begin{abstract}
An $r$-colored permutation is a permutation whose elements have been assigned colors from the set 
$[r]=\{1,2,\dotsc,r\}$ (some colors may be unused), and similarly an $r$-colored rook placement in a Ferrers board
is a rook placement where the rooks are colored with elements of~$[r]$. A colored rook placement then contains a 
colored permutation $\pi$ if it contains a classical occurrence of $\pi$ in which the color of each rook matches
exactly the color of the corresponding entry of~$\pi$.

We say that two $r$-colored permutations $\sigma$ and $\pi$ are \emph{$r$-colored shape-Wilf-equivalent} if 
for every Ferrers board $F$, the number of its $\sigma$-avoiding $r$-colored rook placements 
is equal to the number of the $\pi$-avoiding ones.

We determine the $r$-colored shape-Wilf equivalence for all individual patterns of size up to three, 
as well as for all pairs of patterns of size two, and all values of~$r$. 

As a by-product, we establish several general properties of the $r$-colored shape-Wilf equivalence. For instance, 
we prove that, for every positive integer $n$, all $r$-colorings of the increasing permutation
$(1,2,\dots,n)$ are $r$-colored shape-Wilf equivalent.
\end{abstract}

\section{Introduction}
In this paper, we study a common generalization of two classical concepts from the permutation pattern area. The first
concept is the shape-Wilf equivalence, which deals with occurrences of permutation patterns in rook placements. The 
second concept is the notion of $r$-colored permutations, and its special $r=2$ case of signed permutations.

Let us begin by introducing some basic terminology. A permutation of size $n$ is a sequence 
$\pi=\pi(1),\pi(2),\dotsc,\pi(n)$ in which each element of the set $[n]=\{1, 2, \dots, n\}$ appears once.
We let $S_n$ denote the set of all permutations of size~$n$. When writing small permutations explicitly, we will omit the punctuation and write, e.g., $2143$ instead of $(2,1,4,3)$.

It will often be convenient to represent a permutation $\pi\in S_n$ by its \emph{permutation matrix}, which is an 
$n\times n$ matrix whose entry in column $i$ and row $j$ is equal to $1$ if $\pi(i)=j$ and 0 otherwise. When 
working with matrices and diagrams, we adopt the Cartesian convention of indexing, where rows are 
numbered from bottom to top, columns are numbered from left to right, and by `cell $(i,j)$' we mean the entry
in column $i$ and row $j$.

A permutation $\pi \in S_n$ \emph{contains} a permutation $\sigma \in S_m$, with $n \geq m$ if there is a subsequence of indices $1 \leq i_1 < \dots < i_m  \leq n$ such that $\pi(i_j) < \pi(i_k)$ if and only if $\sigma(j) < \sigma(k)$ for all $1 \leq j, k \leq m$. We say that the positions $\{i_1,\dotsc,i_m\}$ \emph{induce} a copy of $\sigma$ in~$\pi$. For example, the permutation $12543\in S_5$ contains $132\in S_3$, induced by positions $\{1,3,4\}$.
If $\pi$ does not contain $\sigma$, it \emph{avoids} it.

For a permutation $\pi \in S_k$, we let $S_n (\pi)$ denote the set of all permutations of length $n$ that 
avoid~$\pi$. Permutations $\pi$ and $\sigma$ are said to be \emph{Wilf equivalent} if $|S_n (\pi)| = |S_n(\sigma)|$ 
for every~$n\ge 1$.

A \emph{Ferrers board} of order $n$ is a set $F$ of unit squares (`cells') in an $n \times n$ grid arranged into rows and columns such that the rows are left-justified, the row lengths are weakly decreasing from bottom to top, and the number of rows is $n$. 

\begin{figure}[ht!]
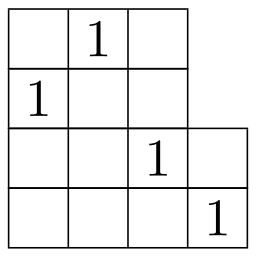

\centering
\young(~1~,1~~,~~1~,~~~1)
\caption{A placement in a Ferrers board. 
Typically, the $0$-entries are omitted to emphasize the positions of the $1$-entries.}
\label{pat:ferrers-board-example}
\end{figure}

A \emph{placement} (also known as \emph{full rook placement} or \emph{transversal}) in a Ferrers board
$F$ is a mapping assigning to each cell of $F$ a value $0$ or $1$, so that each row and each column contains exactly
one entry with value~$1$. See Figure~\ref{pat:ferrers-board-example} for an example.
The entries equal to $1$ are called \emph{rooks}. Notice that a Ferrers board $F$ of order $n$ admits a placement
if and only if its $i$-th row has length at least $n-i+1$. 
Note that permutation matrices correspond to placements in a square $n\times n$ board.

We let $S_F$ denote the set of all placements in a Ferrers board~$F$.

We say that a permutation matrix $B$ (of order $m$) is \emph{contained} in a placement $A$ (of order $n$) 
if there exist two subsets $R := \{r_1 < r_2 < \dots 
< r_m\}$ and $C$ = $\{c_1 < c_2 < \dots < c_m\}$ of the index set $[n]$, such that each of the squares
$(c_i, r_j)$ is within the board of $A$, and $A_{c_i, r_j} = B_{i, j}$ for all $i, j \in [m]$; otherwise we say 
$A$ \emph{avoids}~$B$.  We again say that rows $R$ and columns $C$ \emph{induce} a copy of~$B$ in the placement~$A$.

Notice that a permutation $\pi \in S_n$ contains a permutation $\sigma \in S_m$ if and only if the 
permutation matrix of $\pi$, viewed as a placement in an $n\times n$ board, contains the permutation 
matrix of~$\sigma$. We will often identify a permutation with its matrix, and say that a placement contains 
(or avoids) a permutation $\pi$ if it contains (or avoids) its permutation matrix.

As an example, the permutation $231$ is contained in the placement from Figure~\ref{pat:ferrers-board-example}, 
induced by the rows $R = \{2,3,4\}$ and columns $C = \{1,2,3\}$. On the other hand, the placement avoids the 
permutation $321$. Note that the rows $\{1,2,3\}$ and columns $\{1,3,4\}$ do not induce a copy of $321$, 
since row $3$ does not intersect column $4$ inside the board.

For a permutation $\pi$ and a Ferrers board $F$, we let $S_F(\pi)$ denote the set of all placements in $F$ 
that avoid the permutation $\pi$. Then, we say that two patterns $\pi \in S_n$ and $\sigma \in S_k$ are 
\emph{shape-Wilf-equivalent}, denoted $\pi \swilf \sigma$, if for every Ferrers board $F$, 
$|S_F(\pi)| = |S_F(\sigma)|$. 

Clearly, shape-Wilf equivalence is a refinement of Wilf equivalence, since $|S_n(\pi)|$ is equal to
$|S_F(\pi)|$ where $F$ is the square Ferrers board of order $n$ where each row has~length~$n$.

We now extend these concepts to the colored setting. 

\subsection{Colored shape-Wilf equivalence}

Let $n, r$ be positive integers. An \emph{$r$-colored permutation} of size $n$ is obtained from a permutation 
of size $n$ by assigning to each element a color from the set $[r]$. Formally, we will represent an $r$-colored 
permutation $\pi$ of size $n$ as a pair of maps
\[
\pi_v\colon [n] \to [n], \qquad \pi_c\colon [n] \to [r],
\]
where $\pi_v$ is a permutation of size $n$, and $\pi_c(i)$ is the color assigned to the element~$\pi_v(i)$. 
Note that some of the colors in $[r]$ may remain unused. In particular, any $r$-colored permutation is also
an $(r+1)$-colored permutation. We let $S_n^r$ denote the set of $r$-colored permutations of size~$n$.

For an $r$-colored permutation $\pi=(\pi_v,\pi_c)$, we let $\pi(i)$ denote the value-color pair 
$(\pi_v(i),\pi_c(i))$. 

For convenience in representing small examples, we use the symbols $R$, $B$, $G$ and $Y$
for the colors $1$, $2$, $3$ and $4$, respectively, and refer to them as red, blue, green and yellow.

When writing out small patterns explicitly, we use the shorthand notation ${\pi_v}_{\pi_c}$ for the colored permutation $\pi=(\pi_v,\pi_c)$. 
For example, $\pat{231/RRB}$ denotes the colored permutation $231$, where the two leftmost entries 
(corresponding to values $2$ and $3$) are colored red (i. e., $1$), and the remaining entry is blue.

Containment of colored permutations is defined as in the uncolored case, except that we 
additionally require that the matching entries have the same color. Formally, we say that 
an $r$-colored permutation $\pi=(\pi_v,\pi_c)\in S_n^r$ \emph{contains} an $r$-colored 
permutation $\sigma=(\sigma_v,\sigma_c)\in S_k^r$ if there is a sequence of positions 
$1\le i_1< i_2<\dotsb<i_k\le n$ which induces in $\pi_v$ a copy of $\sigma_v$, and in addition 
satisfies $\pi_c(i_j)=\sigma_c(j)$ for every $j\in[k]$. We let $S_n^r(\sigma)$ denote the set of 
$r$-colored permutations of size $n$ that avoid $\sigma$. We say that two $r$-colored permutations 
$\sigma$ and $\tau$ are
\emph{$r$-colored Wilf-equivalent}, denoted $\sigma\wilf[r]\tau$, if $S_n^r(\sigma)$ and $S_n^r(\tau)$ 
have the same cardinality for every~positive integer~$n$.

An \emph{$r$-colored placement} of a Ferrers board $F$ is obtained from a placement by assigning to each rook 
a color from $[r]$. Formally, it is a mapping that assigns to each cell of $F$ a value $0$ or a pair $(1,c)$ 
for some $c\in[r]$, such that each row and each column contains exactly one entry of the form $(1,c)$.
We let $S^r_F$ denote the set of all $r$-colored placements in a Ferrers board $F$.

We will identify an $r$-colored permutation $\pi\in S_n^r$ with a colored placement in an $n\times n$ board
in the obvious way: if $\pi(i)=(j,c)$, then we place a rook of color $c$ in column $i$ and row~$j$. We refer 
to this placement as the \emph{colored permutation matrix of $\pi$}.

The containment relation of colored permutations and placements refines the uncolored version by requiring that 
the rook that matches an entry of a permutation must have the same color. Formally, given an $r$-colored 
permutation $\pi$ of size $m$ with colored permutation matrix $B$ and a placement $A$ in a Ferrers board $F$ with $n$ columns and rows, 
we say that $\pi$ is \emph{contained} in $A$ if there exist two subsets $R = \{r_1 < r_2 < \dots 
< r_m\}$ and $C = \{c_1 < c_2 < \dots < c_m\}$ of $[n]$, such that each of the squares $(c_i, r_j)$ is 
inside $F$, and $A_{c_i, r_j} = B_{i, j}$ for all $i, j \in [m]$. Otherwise we say that $A$ \emph{avoids}~$B$.

If $P$ is a set of $r$-colored permutations, we say that a placement $A$ \emph{avoids} $P$ if it avoids each permutation in~$P$. 

If $P$ is an $r$-colored permutation or a set of $r$-colored permutations, we write 
 $S^r_F(P)$ for the set of $r$-colored placements in $F$ that avoid~$P$.

Two $r$-colored permutations (or sets of $r$-colored permutations) $P$ and $Q$ are \emph{$r$-colored
shape-Wilf equivalent} if
\[
  |S_F^r(P)|=|S_F^r(Q)|
\]
for every Ferrers board~$F$.  We write this as $P\swilf[r]Q$.

If we identify an uncolored permutation with a $1$-colored permutation whose entries are all colored with
the unique color 1, then the uncolored shape-Wilf equivalence coincides with the 
special case of $r=1$ of the $r$-colored shape-Wilf equivalence. 

For two permutations $\pi\in S_n$ and $\sigma\in S_m$, their \emph{direct sum}, denoted $\pi\oplus\sigma$, is  
the permutation 
\[(\pi(1),\pi(2),\dotsc,\pi(n),n+\sigma(1),n+\sigma(2),\dotsc,n+\sigma(m))\in S_{n+m}.\]
We extend this concept to colored permutations in a straightforward way: if $\pi=(\pi_v,\pi_c)$ 
and $\sigma=(\sigma_v,\sigma_c)$ are colored permutations of size $n$ and $m$, then $\pi\oplus\sigma$ 
is the colored permutation whose values are $\pi_v\oplus\sigma_v$, with the leftmost $n$ elements colored as in $\pi$, and the remaining elements colored as in~$\sigma$. For two sets $P,Q$ of colored permutations, we define
\[
P\oplus Q=\{p\oplus q;\; p\in P, q\in Q\}.
\]

\subsection{Previous results}
While colored shape-Wilf equivalence for a general number of colors does not seem to have been previously studied, it 
is motivated by two concepts that have received considerable attention in the literature, namely the concept of 
shape-Wilf equivalence and the concept of colored (or signed) permutation patterns.

Shape-Wilf equivalence was introduced by Backelin, West, and Xin~\cite{BWX}, who built upon previous 
work of Babson and West~\cite{BabsonWest}. Backelin, West, and Xin~\cite{BWX} showed that the two monotone 
permutations are shape-Wilf equivalent:

\begin{theorem}\label{sw:k1-1k}
    For any positive integer $n$, the permutation $I_n=(1,2, \dots, n)$ is shape-Wilf equivalent to the permutation $D_n=(n, \dots, 2, 1)$.
\end{theorem}

This symmetry between increasing and decreasing patterns has been translated and extended to other 
combinatorial settings, e.g., to more general diagram fillings~\cite{Rubey,Krattenthaler} or to crossings 
and nestings in ordered graphs~\cite{demier}.

Backelin, West and Xin~\cite{BWX} also showed that if a permutation $A$ is shape-Wilf equivalent to a permutation $B$, then for any $C$, 
we also have $A\oplus C$ shape-Wilf equivalent to $B\oplus C$. As we will see in Section~\ref{sec-general}, 
this result generalizes to the colored shape-Wilf setting.

Another result on shape-Wilf equivalence, due to Stankova and West~\cite{StankovaWest2002}, states that the 
two permutations 231 and 312 are shape-Wilf equivalent. At least two other proofs of this result are 
known~\cite{JELINEK2007,BloomSaracino}, and versions of this equivalence have also been considered in other combinatorial settings, such as words~\cite{GKZ} or partial permutations~\cite{CJJK}.

These results allow a complete classification of permutations of size three with respect to shape-Wilf equivalence, 
with the three equivalence classes being $\{123,321,213\}$, $\{231,312\}$ and $\{132\}$~\cite{JELINEK2007}. 
We note that Bloom and Elizalde~\cite{BE} have also identified the shape-Wilf classes of pairs of patterns of size~3.

Colored permutations and their containment and Wilf-classification go back at least to the work of 
Mansour~\cite{mansour}, who performed a systematic study of $r$-colored (non-shape) Wilf-equivalence 
among colored patterns of size~2, as well as their two-element subsets. Other authors considered 
specific classes of colored pattern sets~\cite{mansour2,PRW,Egge}.

In the area of pattern avoidance in colored permutations, by far the most frequently studied case
is the case of two colors, which is usually presented as \emph{signed permutations} where instead of 
distinguishing elements by two colors, they are distinguished by the two signs $+$ and~$-$. From the 
point of view of Wilf classification, the two settings are however equivalent. 

Pattern-avoidance in signed permutations has been considered by Simion~\cite{Simion}, and by Billey and
Lam~\cite{BilleyLam}. Mansour and West~\cite{MansourWest} have systematically classified all
the subsets of~$S_2^2$ with respect to 2-colored Wilf equivalence. Dukes, 
Mansour and Reifegerste~\cite{DMR2008} and Dukes, Jelínek, 
Mansour and Reifegerste~\cite{DJMR2009} 
have classified signed permutations of size up to 5 into Wilf classes, and also addressed 
shape-Wilf equivalence. Among other results, Dukes et al.~\cite[Corollary 2.3]{DMR2008} showed that 
all the $2^n$ possible assignments of signs to the increasing permutation $I_n=(1,2,3,\dotsc,n)$ 
yield Wilf-equivalent patterns, and this was later extended to shape-Wilf 
equivalence~\cite[Theorem 3.1]{DJMR2009}. In our paper, we further extend this result to any number 
of colors; see Proposition~\ref{prop:1n}.

Other authors have applied the coloring approach to settings beyond classical permutation containment, 
like set partitions~\cite{GoytPudwell,GoytPudwell2}, circular permutations~\cite{GrayLanningWang}, consecutive 
patterns~\cite{DuaneRemmel} or matchings~\cite{ChenGuo}.

Recently, Blanco and Skora~\cite{BlancoSkora} have studied the possible growth function 
of hereditary sets of colored permutations, focusing on classes of polynomial growth. 

\subsection*{Our results}

The main focus of this paper is on $r$-colored shape-Wilf equivalence, where we obtain a complete 
classification of $r$-colored shape-Wilf equivalence for colored patterns of size up to 3, and for
pairs of patterns of size 2. Our classification is valid for any value of $r$; in fact, a general result
(Lemma~\ref{lem:color-extension}) shows that the $r$-shape-Wilf equivalence does not depend 
on the choice of palette size~$r$, provided the palette $[r]$ contains all the colors that appear in the 
patterns considered. Therefore, we may limit ourselves to the situation when the palette size used 
for the rook placements is just large enough to contain all the colors in the avoided patterns,
since a smaller palette size would mean that all the colored rook placements avoid a pattern we consider, 
while a larger palette size will display no new behavior in terms of colored shape-Wilf equivalence. 

As a corollary to our results on colored shape-Wilf equivalence, we also obtain a full classification 
of colored Wilf equivalence for the same pattern sets. It will in fact turn out that for colored 
patterns of size at most 3, as well as for pairs of colored 
patterns of size 2, we may recover all the $r$-colored Wilf classes by simply taking the corresponding 
$r$-colored shape-Wilf classes and merging those that can be related by the standard symmetries of the square.
These results extend the previous classifications of signed patterns of size 3~\cite{DMR2008} to patterns with
any number of colors, while the Wilf classes
of pairs of colored patterns have been already determined by Mansour~\cite{mansour} by other methods.

In Section~\ref{sec-general}, we derive general results and develop useful tools for proving
colored shape-Wilf equivalences. Next, in
Section~\ref{sec:singleton}, we classify singleton
$r$-colored patterns of length at most three under
$r$-colored shape-Wilf equivalence, for every integer
$r\geq 2$. In Section~\ref{sec:sets-of-patterns}, we perform a similar classification for
sets containing two colored patterns of size~2.


\section{General results}\label{sec-general}
Before we deal with specific instances of colored Wilf equivalence, we establish some general facts that
will be useful throughout the rest of this paper.

\subsection{Irrelevant colors}
We first show that shape-Wilf equivalence of colored patterns is unaffected by adding colors that do not appear in any 
of the patterns.
\begin{lemma}
\label{lem:color-extension}
    Let $P$ and $Q$ be sets of $r$-colored patterns. Let $k \geq r$ be a positive integer. Then 
    \[P \swilf[r] Q \iff P \swilf [k] Q, \qquad\text{and}\qquad P \wilf[r] Q \iff P \wilf [k] Q.\]
\end{lemma}
\begin{proof}
We only deal with the first equivalence, as the proof of the second one is analogous.

$\Longrightarrow\colon$ Suppose $P \swilf[r] Q$. Let $F$ be an arbitrary Ferrers board and let $N 
\in S_F^{k}(P)$ be a placement avoiding~$P$. Fix all the rows and columns of elements in $N$ that 
have any color from $\{r + 1, \dots, k\}$. Then remove their rows and columns; we obtain a 
sub-placement $N'$ that induces a board $F'$. $N'$ avoids $P$, otherwise $N$ would contain 
$P$. Moreover, $N'$ contains only colors from $\{1, \dots, r\}$. Then since $P \swilf[r] Q$, there is 
a bijection $f_{F'}:S^{r}_{F'}(P) \to S^{r}_{F'}(Q)$. We can apply it to~$N'$. Now, we replace $N'$ 
in $N$ with $f_{F'}(N')$. The new placement avoids~$Q$. The construction is 
invertible since the rooks outside $N'$ remain unchanged; for the inverse we use $f_{F'}^{-1}$.

$\Longleftarrow\colon$ Suppose $P \swilf[k] Q$. We will prove that for every Ferrers board $F$ 
$|S^r_F(P)| = |S^r_F(Q)|$ by induction on the order of~$F$. For the empty board, both counts are~$1$. 
The case $k = r$ is immediate, so we can assume $k > r$. Let us say that a set of cells $E\subseteq F$ 
is \emph{admissible} if $F$ has a rook placement in which all the cells of $E$ contain rooks.

For an admissible set $E$, let $F/E$ denote the Ferrers board obtained by removing all 
the rows and columns from $F$ that contain a cell in $E$, and by closing the gaps. Note that admissibility 
of $E$ guarantees that if in the above deletion process we remove all the cells in a column $i$ of $F$ because their
rows contain a cell from $E$, then column $i$ itself must contain a cell from $E$ and will also be deleted. 
In other words, $F/E$ does not contain `columns of zero length' that evaded deletion. Symmetrically, $F/E$ 
does not contain rows of zero length either.

We may group the 
placements in $S^k_F(P)$ by the set $E$ of rooks that are colored by the `extra' colors in 
$[k]\setminus[r]$. Note that these rooks can never take part in any copy of a forbidden pattern. 

This shows that
\[|S^k_F(P)| = \sum_{E} |k-r|^{|E|}|S^{r}_{F/E}(P)|,\]
where the summation ranges over all the admissible sets $E\subseteq F$.

When $E = \emptyset$, in the sum we get the term $|S^r_F(P)|$. 

Similarly, we can count the $Q$-avoiding placements, again summing over all admissible $E\subseteq F$, as
\[|S^k_F(Q)| = \sum_{E} |k-r|^{|E|}|S^{r}_{F/E}(Q)|.\]
Again, when $E = \emptyset$ in the sum we get the term $|S^r_F(Q)|$.

By induction hypothesis, $|S^r_{F/E}(P)| = |S^r_{F/E}(Q)|$ for $E \neq \emptyset$. Thus 
\begin{align*}
0 &= |S^k_F(P)| - |S^k_F(Q)|\\& = \sum_{E} |k-r|^{|E|}(|S^{r}_{F/E}(P)| - |S^{r}_{F/E}(Q)|)\\ &= |S^r_F(P)| - |S^r_F(Q)|. \qedhere
\end{align*}
\end{proof}
In view of Lemma~\ref{lem:color-extension}, when studying shape-Wilf equivalence of colored patterns, we will rarely 
care about the precise size of the color palette $[r]$, as long as it is large enough to contain all the colors used by 
the patterns. Henceforth, for sets of patterns $P$ and $Q$ we will write $P\swilf[*]Q$ as shorthand for $P\swilf[r] Q$
where $r$ is the smallest integer such that all the colors used by the patterns in $P\cup Q$ are in~$[r]$. 
Lemma~\ref{lem:color-extension} then implies that $\swilf[*]$ is an equivalence relation on sets of colored patterns.
We refer to this relation simply as the \emph{colored shape-Wilf equivalence}. We apply the same convention to 
colored Wilf equivalence, which will be denoted $P\wilf[*] Q$.

\subsection{Easy symmetries}

Let $c\colon[r]\to[r]$ be a bijection. For an $r$-colored permutation $\pi=(\pi_v,\pi_c)$, we let $c(\pi)$ denote the 
$r$-colored permutation obtained from $\pi$ by changing each color $\alpha\in[r]$ into $c(\alpha)$. Formally speaking, 
$c(\pi)$ is the colored permutation $\sigma=(\sigma_v,\sigma_c)$, where $\sigma_v=\pi_v$ and $\sigma_c(i)=c(\pi_c(i))$ 
for every $i\in[n]$. Similarly, for a set of $r$-colored patterns $P$, we let $c(P)$ denote the set $\{c(\pi);\; \pi\in P\}$.

We say that two $r$-colored patterns $\sigma$ and $\pi$ are \emph{color-isomorphic} if there is a bijection 
$c\colon[r]\to[r]$ such that $\sigma=c(\pi)$; similarly, two sets $P$ and $Q$ of $r$-colored patterns are 
color-isomorphic if $Q=c(P)$ for a bijection $c\colon[r]\to[r]$. The following fact has already been observed by Mansour~\cite[Remark 2.4]{mansour} in the context of (non-shape) colored Wilf equivalence.

\begin{observation}\label{color-isomorphism}
Two color-isomorphic patterns, or sets of patterns, are colored shape-Wilf equivalent. 
\end{observation}

Beyond color-isomorphism, another useful symmetry operation is the transpose, i.e., the reflection along the $x=y$ axis. Note that the transpose of a Ferrers board is again a Ferrers board. For a Ferrers board, placement, or pattern $X$, we let $X^t$ denote its transpose. If $P$ is a set of patterns, we let $P^t$ denote the set $\{\pi^t;\; \pi\in P\}$.

\begin{observation}\label{obs-transpose}
If $X$ and $Y$ are colored shape-Wilf equivalent patterns, or sets of patterns, then $X^t$ and $Y^t$ are colored 
shape-Wilf equivalent as well.
\end{observation}

We remark that a colored pattern is not necessarily colored shape-Wilf equi\-valent to its transpose. For instance, 
\pat{321/RRB} is not colored shape-Wilf equivalent to \pat{321/BRR}, as we will see in Section~\ref{sec:singleton}.

\subsection{Direct sums}
Many known examples of (uncolored) shape-Wilf equivalence are based on a general argument due to Backelin, West, and 
Xin~\cite[Proposition~2.3]{BWX}, later explicitly extended to sets of patterns by Burstein et al.~\cite{BHKZ}, which
allows one to bootstrap a single pair of equivalent patterns (or sets of patterns) 
into an infinite family of equivalent pairs using direct sums. The argument can easily be adapted to the realm of 
colored patterns, as we now show.

Let $T$ be a colored placement on a Ferrers board $F$, and let $Q$ be a nonempty set of colored patterns.
We construct the \emph{$Q$-marking} of $F$ as follows.
\begin{enumerate}
  \item Mark a cell $(i,j)$ if the subboard strictly northeast of $(i,j)$
        contains an occurrence of a member of $Q$.
  \item Whenever a rook of $T$ lies in an unmarked cell, unmark its entire
        row and column.
  \item Delete the unmarked rows and columns and close the resulting gaps.
        Denote the board formed by the marked cells by $F_Q(T)$ and the
        restriction of $T$ to its marked cells by $T_Q$.
\end{enumerate}

The board $F_Q(T)$ is a Ferrers board, and $T_Q$ is a full colored placement
on it.  The reason for using the northeast subboard is that, with our matrix
convention, the second block of a direct sum lies northeast of the first
block.  

Consequently, we directly observe that
\begin{equation}\label{eq:coloring-criterion}
T\text{ contains a member of }P\oplus Q
\quad\Longleftrightarrow\quad
T_Q\text{ contains a member of }P.
\end{equation}

\begin{figure}[htbp]
\centering
\begin{tikzpicture}[
scale=0.7, transform shape,
x=0.54cm,y=0.54cm,
  boardcell/.style={draw=black!65,line width=0.35pt},
  rook/.style={circle,inner sep=0pt,minimum size=11pt,
    font=\scriptsize\bfseries,text=white},
  panel/.style={font=\small\bfseries},
  index/.style={font=\scriptsize,text=black!70}]
\node[font=\small] at (13.1,10.1)
  {$Q=\patset{21/RB},\quad F=(7,7,7,7,6,5,4),
    \quad T=\pat{1576243/RRRBBRB}$};
\node[panel] at (3.5,8.7) {1. Mark cells};
\node[panel] at (13.5,8.7) {2. Unmark rows and columns};
\node[panel] at (23.5,8.7) {3. Close gaps};

\foreach \stage/\offset in {1/0,2/10} {
  \begin{scope}[shift={(\offset,0)}]
    \foreach \height [count=\col] in {7,7,7,7,6,5,4} {
      \foreach \row in {1,...,\height} {
        \pgfmathtruncatemacro{\ismarked}
          {(\col<=2 && \row<=5) || (\col<=5 && \row<=2)}
        \pgfmathtruncatemacro{\iskept}
          {(\col==1 || \col==2 || \col==5)
            && (\row==1 || \row==2 || \row==5)}
        \def\cellfill{white}
        \ifnum\ismarked=1\relax
          \ifnum\stage=1\relax
            \def\cellfill{marked}
          \else
            \ifnum\iskept=1\relax
              \def\cellfill{marked}
            \fi
          \fi
        \fi
        \path[fill=\cellfill]
          ({\col-1},{\row-1}) rectangle (\col,\row);
        \ifnum\stage=2\relax
          \ifnum\ismarked=1\relax
            \ifnum\iskept=0\relax
              \path[pattern=north east lines,pattern color=black!30]
                ({\col-1},{\row-1}) rectangle (\col,\row);
            \fi
          \fi
        \fi
        \draw[boardcell]
          ({\col-1},{\row-1}) rectangle (\col,\row);
      }
    }
    \foreach \index in {1,...,7} {
      \node[index] at ({\index-0.5},-0.5) {\index};
      \node[index] at (-0.45,{\index-0.5}) {\index};
    }
    \foreach \row/\shade/\letter [count=\col] in
      {1/rookR/R,5/rookR/R,7/rookR/R,6/rookB/B,
       2/rookB/B,4/rookR/R,3/rookB/B} {
      \pgfmathtruncatemacro{\rookkept}
        {\col==1 || \col==2 || \col==5}
      \ifnum\stage=2\relax
        \ifnum\rookkept=0\relax
          \node[rook,fill=\shade!25,text=black!65]
            at ({\col-0.5},{\row-0.5}) {\letter};
        \else
          \node[rook,fill=\shade]
            at ({\col-0.5},{\row-0.5}) {\letter};
        \fi
      \else
        \node[rook,fill=\shade]
          at ({\col-0.5},{\row-0.5}) {\letter};
      \fi
    }
    \ifnum\stage=2\relax
      \foreach \index in {3,4,6,7} {
        \draw[black!55,line width=0.6pt]
          ({\index-0.7},-0.7) -- ({\index-0.3},-0.3);
        \draw[black!55,line width=0.6pt]
          (-0.65,{\index-0.7}) -- (-0.25,{\index-0.3});
      }
    \fi
  \end{scope}
}
\draw[dashed,line width=0.8pt]
  (2.06,5.06) rectangle (3.94,6.94);
\draw[dashed,line width=0.8pt]
  (5.06,2.06) rectangle (6.94,3.94);
\draw[-{Latex[length=2mm]},black!60,line width=0.7pt]
  (7.6,3.5) -- (9.1,3.5);
\draw[-{Latex[length=2mm]},black!60,line width=0.7pt]
  (17.6,3.5) -- (21.0,3.5);
\node[font=\small] at (23.5,6.2) {$F_Q(T)=(3,3,2)$};
\begin{scope}[shift={(22,2)}]
  \foreach \height [count=\col] in {3,3,2} {
    \foreach \row in {1,...,\height} {
      \path[fill=marked,draw=black!65,line width=0.35pt]
        ({\col-1},{\row-1}) rectangle (\col,\row);
    }
  }
  \foreach \index in {1,2,3} {
    \node[index] at ({\index-0.5},-0.5) {\index};
    \node[index] at (-0.45,{\index-0.5}) {\index};
  }
  \foreach \row/\shade/\letter [count=\col] in
    {1/rookR/R,3/rookR/R,2/rookB/B} {
    \node[rook,fill=\shade] at ({\col-0.5},{\row-0.5}) {\letter};
  }
\end{scope}
\node[font=\small] at (23.5,0.5) {$T_Q=\pat{132/RRB}$};
\node[font=\scriptsize,align=center] at (3.5,-1.8)
  {Dashed rectangles show\\occurrences of $Q$.};
\node[font=\scriptsize,align=center] at (13.5,-1.8)
  {Keep rows and columns\\$1,2,5$.};
\path[fill=marked,draw=black!65]
  (4.2,-3.4) rectangle (4.7,-2.9);
\node[anchor=west,font=\scriptsize] at (4.95,-3.15) {marked cell};
\path[pattern=north east lines,pattern color=black!30,draw=black!65]
  (10.0,-3.4) rectangle (10.5,-2.9);
\node[anchor=west,font=\scriptsize] at (10.75,-3.15)
  {mark removed in step 2};
\node[rook,fill=rookR!25,text=black!65] at (19.1,-3.15) {R};
\node[anchor=west,font=\scriptsize] at (19.55,-3.15)
  {rook to be deleted};
\end{tikzpicture}
\caption{The $Q$-marking of a colored placement.
Board tuples give row lengths from bottom to top. Gray cells are marked.
The initially unmarked rooks lie in rows $3,4,6,7$ and columns $3,4,6,7$.
Clearing these rows and columns and closing the gaps leaves
$F_Q(T)=(3,3,2)$ with $T_Q=\pat{132/RRB}$.}
\label{fig:q-coloring}
\end{figure}
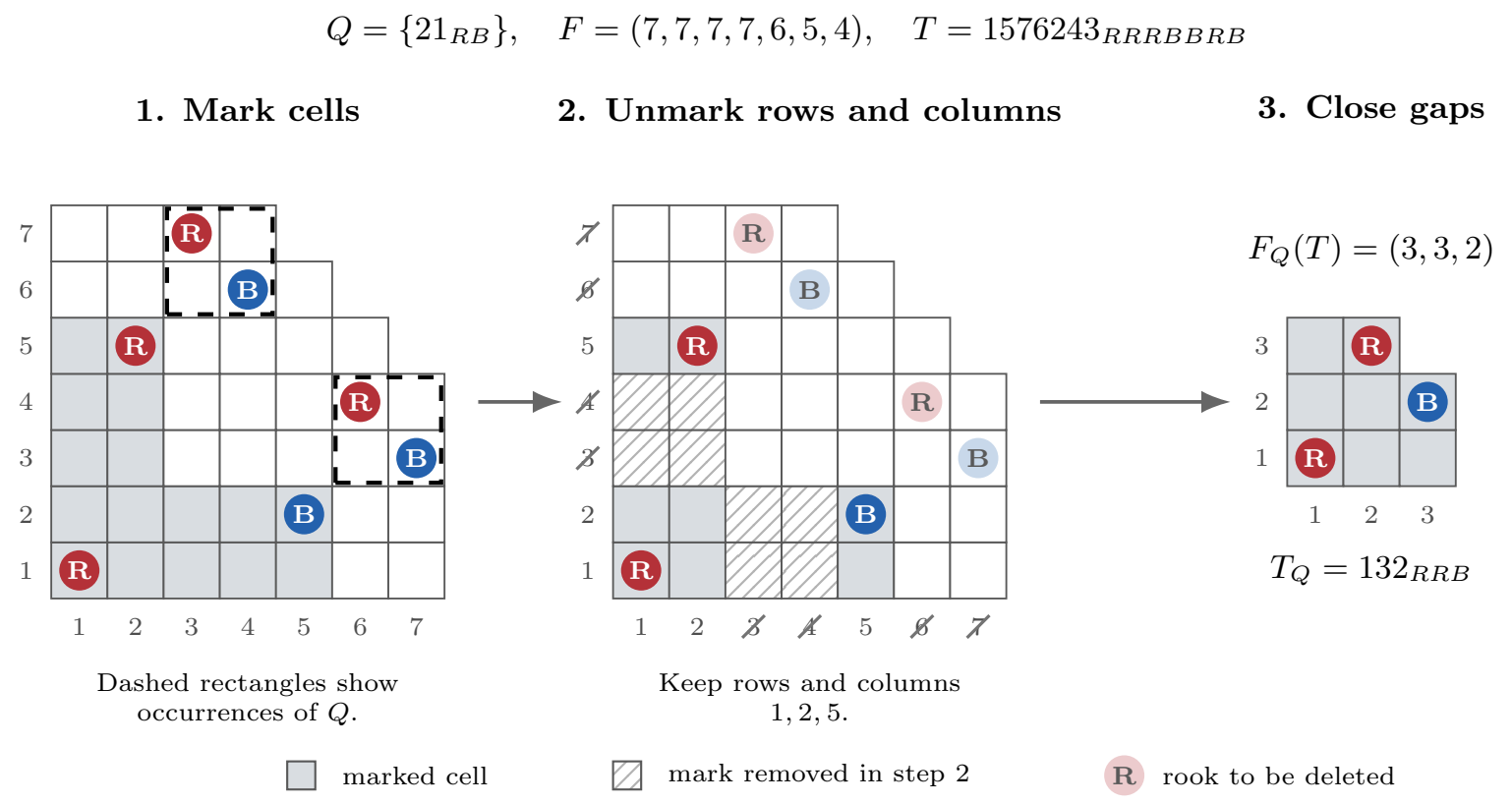

Note that replacing $T_Q$ by another colored placement on the same marked
board does not change the $Q$-marking; indeed, we may easily see that any marked cell has 
its marking witnessed by an occurrence
of a pattern from $Q$ that only uses unmarked cells. Thus, by modifying the placement in the marked cells, we can
neither create a new marked cell, nor eliminate every witness of an existing marked cell.

\begin{lemma}\label{lem:direct-sum-extension} Let $P$, $P'$ and $Q$ be nonempty sets of $r$-colored patterns.
If $P\swilf[*]P'$, then $P\oplus Q\swilf[*]P'\oplus Q$.
\end{lemma}

\begin{proof}
For every Ferrers board $G$, fix a bijection
\[
  \varphi_G\colon S_G^r(P)\longrightarrow S_G^r(P').
\]
Let $T\in S_F^r(P\oplus Q)$.  Construct $F_Q(T)$ and $T_Q$ as above.  By
\eqref{eq:coloring-criterion}, $T_Q$ avoids~$P$.  Replace $T_Q$ by
$\varphi_{F_Q(T)}(T_Q)$ and leave all unmarked rooks unchanged.  The set of $Q$-marked cells
does not change, so another application of
\eqref{eq:coloring-criterion} shows that the resulting placement avoids
$P'\oplus Q$.

Using $\varphi_G^{-1}$ on every marked board reverses the construction.
Hence the construction is a bijection for every $F$.
\end{proof}

We remark that the order of the two summands in the previous lemma cannot, in general, be reversed:
$P\swilf[*]P'$ does not imply $Q\oplus P\swilf[*]Q\oplus P'$, even in the classical uncolored case of 
$r=1$. For instance, the permutations $12$ and $21$ are known to be shape-Wilf equivalent, while 
$1\oplus12 =123$ and $1\oplus21=132$ are known not to be.

We may, however, show that colored shape-Wilf equivalence of direct sums is preserved by color 
isomorphism of the second summands.

\begin{lemma}
\label{lem:tail-recoloring}
Let $P$, $Q$ and $Q'$ be sets of colored patterns, and suppose that $Q$ is color-isomorphic to $Q'$. Then
\[
  P\oplus Q\swilf[*]P\oplus Q'.
\]
\end{lemma}

\begin{proof}
Suppose that $P$, $Q$ and $Q'$ are $r$-colored, and
let $c\colon[r]\to[r]$ be the bijection such that $Q'=c(Q)$.
We have \[
    P \oplus Q \swilf[*] c(P \oplus Q) = c(P)\oplus Q' \swilf[*] P \oplus Q',
\]
where the first equivalence is due to Observation~\ref{color-isomorphism}, and the second equivalence is due 
to Lemma~\ref{lem:direct-sum-extension} combined with Observation~\ref{color-isomorphism}. \end{proof}

Combining Lemma~\ref{lem:direct-sum-extension} and Lemma~\ref{lem:tail-recoloring}, we reach the following conclusion.

\begin{corollary}\label{cor:two-changes}
Let $P, P'$, $Q$ and $Q'$ be sets of colored patterns such that $P \swilf[*] P'$ and $Q$ is color-isomorphic to $Q'$. 
Then $P \oplus Q \swilf[*] P' \oplus Q'$.
\end{corollary}

\subsection{Colorings of the identity permutation}

Let $\iota_n=(1,2,\dotsc,n)$ be the identity permutation of length~$n$. We will show that all the colorings of $
\iota_n$ are colored shape-Wilf equivalent. This generalizes a result of Dukes, Mansour and
Reifegerste~\cite[Corollary 2.3]{DMR2008}, who established the special case of two colors (i. e., signed 
permutations) and (non-shape) Wilf equivalence, which was then extended to shape-Wilf equivalence by 
Dukes, Jelínek, Mansour and Reifegerste~\cite[Theorem 3.1]{DJMR2009}. For more than two colors,
Mansour~\cite{mansour} showed that all colorings of $\iota_2=12$ are Wilf equivalent. 

\begin{proposition}
\label{prop:1n}
For any positive integer $n$, any two colorings of the increasing permutation $\iota_n$ are colored shape-Wilf 
equivalent.
\end{proposition}

\begin{proof}
Let $r$ be a positive integer. Let us proceed by induction on $n$. The base case $n = 1$ trivially holds.

For the inductive step, let $n>1$ and let $\pi$ and $\sigma$ be $r$-colored permutations of size $n$ such that 
$\pi_v = \sigma_v= \iota_n$. Let $\alpha$ be the color of the rightmost entry in $\pi$ and $\beta$ the color of the 
rightmost entry in $\sigma$. In particular we can write $\pi$ as $\pi'\oplus 1_\alpha$ and $\sigma$ as $\sigma'\oplus 
1_\beta$, where $\pi'$ and $\sigma'$ are two colorings of $\iota_{n-1}$. 

Clearly $1_\alpha$ is color-isomorphic to $1_\beta$, and by induction, $\pi'\swilf[*] \sigma'$. By 
Corollary~\ref{cor:two-changes}, we conclude that $\pi\swilf[*]\sigma$, as claimed.
\end{proof}

\subsection{Disjoint palettes}

\begin{lemma}\label{lem:pattern-addition}
Let $P$, $P'$ and $Q$ be sets of colored patterns, where the set of colors appearing in the patterns of $P\cup P'$ is disjoint from the set of colors appearing in~$Q$. If $P\swilf[*] P'$, then $P\cup Q\swilf[*] P'\cup Q$.
\end{lemma}
\begin{proof}
Suppose that the patterns in $P\cup P'$ use colors from the set $[t]$, while those in $Q$ use colors from the set 
$[r]\setminus[t]$, for some $1\le t<r$. Fix a Ferrers board~$F$. Our aim is to find a bijection between the 
$r$-colored placements in $F$ that avoid $P\cup Q$ and those that avoid $P'\cup Q$.

Suppose we are given an $r$-colored placement $T$ of a board $F$ avoiding $P \cup Q$. 
Let $A$ be the set of all the rooks of $T$ that are colored by colors from the set 
$[t]$, and let $B$ be the set of all the remaining rooks. Let $F_A$ be the sub-diagram 
of $F$ induced by the rows and columns containing the rooks of $A$, and define $F_B$ analogously. 
After closing the gaps, $F_A$ and $F_B$ are Ferrers boards. Let $T_A$ and $T_B$ denote the
restrictions of $T$ to $F_A$ and $F_B$, respectively. These restrictions
are full colored placements on the corresponding boards.
Then we have \[
  T\text{ avoids }P\cup Q
  \quad\Longleftrightarrow\quad
  T_A\text{ avoids }P
  \text{ and }
  T_B\text{ avoids }Q.
\]

Since $P\swilf[t]P'$, we may choose a bijection that transforms $T_A$ into
a $t$-colored placement on $F_A$ avoiding~$P'$. Leaving $T_B$ unchanged, we obtain an
$r$-colored placement on $F$ avoiding $P'\cup Q$. Applying the inverse
bijection to the restriction on $F_A$ reverses the construction. Hence this
is a bijection, and therefore
\[
  P\cup Q\swilf[*]P'\cup Q.\qedhere
\]
\end{proof}

\subsection{Extreme elements in placements}\label{ssec-extreme}
We now derive an auxiliary result on equidistribution of certain statistics in uncolored placements.

Consider an (uncolored) placement $T$ of a Ferrers board $F$. We say that two rooks 
$r_1$ and $r_2$ of $T$ are \emph{compatible} if $F$ contains a rectangular 
subboard that contains both $r_1$ and~$r_2$. In other words, $r_1$ and $r_2$ are compatible if they 
form an occurrence of the pattern 12 or 21. 

Suppose that $r_1$ and $r_2$ are two distinct compatible rooks, with $r_1$ in column $i_1$ 
and row $j_1$, and $r_2$ in column $i_2$ and row~$j_2$. We say that $r_1$ is to the North-West of $r_2$ (and 
$r_2$ is to the South-East of $r_1$) if $i_1<i_2$ and $j_1>j_2$. South-West and North-East are defined 
analogously. 

We say that a rook $r$ is \emph{NW-extreme} in $T$ if $T$ does not contain any other 
rook compatible with $r$ to the North-West of~$r$ (informally speaking, $r$ never plays the role of `1' in an 
occurrence of 21). Again, NE-extreme, SW-extreme and SE-extreme 
rooks are defined analogously. We let $\NE(T)$ denote the number of NE-extreme rooks in $T$, and 
similarly for $\NW(T)$, $\SE(T)$ and $\SW(T)$.

We will need the following equidistribution result, which is proved by adapting an argument of 
Poznanovi\'c~\cite[Theorem 2.2]{Poznanovic} on ordered matchings, via the well-known correspondence 
between rook placements of Ferrers shapes and ordered matchings (see, e. g., de Mier~\cite{demier}).

\begin{theorem}\label{thm-equi}
For any Ferrers shape $F$, the three statistics $\NW(\cdot)$, $\SW(\cdot)$ and $\SE(\cdot)$ have the 
same distribution on the set of placements of~$F$; in other words, for any integer $k$, the 
number of placements of $F$ with $k$ NW-extreme rooks is the same as the number of those with $k$ 
SW-extreme rooks and those with $k$ SE-extreme rooks. 
\end{theorem}
\begin{proof}
Let $F$ be a Ferrers diagram with $n$ rows and $n$ columns, and let $\lambda_j$ be the length of 
the $j$-th row of $F$, so that we have 
\[
 n=\lambda_1\ge \lambda_2\ge\dotsb\ge \lambda_n.
\]
We may assume that $\lambda_j\ge n-j+1$, otherwise $F$ would have no placement. 

Let $a_{pq}$ be the number of placements of $F$ that have $p$ SW-extreme rooks and $q$ SE-extreme 
rooks. We will show that $a_{pq}=a_{qp}$, implying that $\SW(\cdot)$ and $\SE(\cdot)$ have symmetric 
joint distribution. A transposed version of the same argument then yields that $\SW(\cdot)$ and 
$\NW(\cdot)$ also have symmetric joint distribution, and the theorem follows.

To establish the symmetry of $\SW(\cdot)$ and $\SE(\cdot)$, consider the generating polynomial 
\[
 P(x,y)=\sum_{p,q\ge 0} a_{p,q}x^p y^q,
\]
with the aim of showing that $P(x,y)=P(y,x)$. We will derive an explicit product formula for 
$P(x,y)$ displaying the required symmetry directly.

Let us generate a placement of $F$ by processing rows from top to bottom, and placing a rook in 
each row $j$ to a column distinct from all the columns that already contain the rooks in rows above 
row~$j$. In particular, in row $j$ we may choose among $\lambda_j-n+j$ possible positions. The 
key observation is that the rook we place in row $j$ will end up being SW-extreme in the completed placement 
if and only if we choose the leftmost available position, and it will end up being SE-extreme if and only if we choose 
the rightmost available position. Note that when $\lambda_j=n-j+1$, we only have one available position in 
row $j$, which will yield a rook that is both SW-extreme and SE-extreme. It follows that
\[
P(x,y)=\prod_{j=1}^n R_j(x,y), 
\]
where 
\[
 R_j(x,y)=\begin{cases}
           x+y+\lambda_j-n+j-2 &\text{if $\lambda_j>n-j+1$}\\
           xy &\text{if $\lambda_j=n-j+1$}.
          \end{cases}
\]
In any case, $R_j(x,y)$ is equal to $R_j(y,x)$, and the theorem follows.
\end{proof}

As mentioned in the proof above, it actually shows that $\SW(\cdot)$ and $\SE(\cdot)$ have symmetric joint
distribution on placements of a given Ferrers board, rather than merely having identical distributions. It 
even shows that the same symmetry applies to the joint distribution of the set-valued statistics
$\rowSE(\cdot)$ and $\rowSW(\cdot)$ defined as the sets of indices of the rows containing $\SE$-extreme and 
$\SW$-extreme rooks, respectively. Moreover, this symmetry of $\rowSE(\cdot)$ and $\rowSW(\cdot)$ holds,
with the same proof, 
to placements in (bottom-aligned) stack polyominoes, which are a generalization of Ferrers 
boards that have been previously studied in the context of rook placements~\cite{GuoPoznanovic,Jonsson}. 
We, however, do not need this level of generality in our paper.


\section{Avoidance of single patterns of size up to 3}\label{sec:singleton}

We will now present the complete list of colored shape-Wilf classes of individual colored permutations of 
size up to 3. We will also point out that by collapsing patterns using the dihedral symmetries of the square, 
our classification recovers the known Wilf-equivalence classification obtained for colored patterns of size 2 by 
Mansour~\cite{mansour} and for signed patterns of size 3 by Dukes, Mansour and Reifegerste~\cite{DMR2008}. 

Since color-isomorphic patterns are colored shape-Wilf equivalent, we only need to consider 
one representative from each color-isomorphism class. We will say that an $r$-colored permutation $\pi$ is 
\emph{canonical} if for every color $\alpha>1$ that appears in $\pi$, the color $\alpha-1$ appears as well and the 
leftmost element colored by $\alpha-1$ is to the left of the leftmost element colored by~$\alpha$.

Clearly, every colored permutation is color-isomorphic to a unique canonical one, so we only need to consider 
canonical patterns in our classification. We will, for convenience, refer to colors by these codes: $R=1$, $B=2$, $G=3$ and $Y=4$, referring to them as red, blue, green and yellow. There is only one canonical colored permutation of size 1, namely $1_R$, and nothing more to classify. For sizes two and three, Table~\ref{tab:singleton-classes} lists all the canonical patterns and their colored shape-Wilf equivalence classes.

\begin{table}[ht]
\centering
\small
\renewcommand{\arraystretch}{1.25}
\begin{tabular}{@{}cl}
\toprule
Class & Patterns \\
I
& \pat{12/RR}, \pat{12/RB}, \pat{21/RR}\\

II
& $\pat{21/RB}$ \\

\midrule
III
& \pat{123/RRR}, \pat{123/RRB}, \pat{123/RBR}, \pat{123/RBB}, \pat{123/RBG}, \pat{213/RRR}, \pat{213/RRB}, \pat{321/RRR}\\

IV
& \pat{321/RRB} \\

V
& \pat{321/RBR}\\

VI
& \pat{321/RBB}\\

VII
& \pat{213/RBR}, \pat{213/RBB}, \pat{213/RBG}\\

VIII
& \pat{231/RRR}, \pat{312/RRR}\\

IX
& \pat{132/RRR}, \pat{132/RBB}\\

X
& \pat{132/RRB}, \pat{132/RBR}, \pat{132/RBG} \\

XI
& \pat{231/RRB} \\

XII
& \pat{231/RBR} \\

XIII
& \pat{231/RBB} \\

XIV
& \pat{312/RRB} \\

XV
& \pat{312/RBR}\\

XVI
& \pat{312/RBB} \\

XVII
& \pat{321/RBG}\\

XVIII
& \pat{231/RBG}\\

XIX
& \pat{312/RBG} \\
\bottomrule
\end{tabular}
\caption{The $3$-colored shape-Wilf equivalence classes of
canonical singleton patterns of lengths two and three.
Each row describes one equivalence class.}
\label{tab:singleton-classes}
\end{table}

\begin{theorem}\label{thm-singleton}
The nineteen rows of Table~\ref{tab:singleton-classes} correspond to the colored shape-Wilf equivalence classes of canonical colored permutations of size 2 and 3.
\end{theorem}
\begin{proof}
Let us begin by arguing that patterns listed on the same row are indeed equivalent. 

\emph{Class I.} The equivalence of \pat{12/RR} and \pat{21/RR} follows from the shape-Wilf equivalence of 12 and 21, and the equivalence of \pat{12/RR} and \pat{12/RB} follows from Proposition~\ref{prop:1n}. 

\emph{Class III.} The first five patterns are equivalent by Proposition~\ref{prop:1n}, the three monochromatic ones are equivalent because of the shape-Wilf equivalence of the uncolored permutations 123, 213 and 321, and finally, \pat{213/RRR} is equivalent to \pat{213/RRB} by Lemma~\ref{lem:tail-recoloring}. 

\emph{Class VII.} All three patterns are equivalent by Lemma~\ref{lem:tail-recoloring}.

\emph{Class VIII.} The two patterns are equivalent, because the corresponding uncolored permutations are shape-Wilf equivalent.

\emph{Class IX.} The two patterns are equivalent by Lemma~\ref{lem:tail-recoloring}.

\emph{Class X.} All three patterns are equivalent by Lemma~\ref{lem:tail-recoloring}.

It remains to argue that patterns in distinct rows of the table are not equivalent. For this, it suffices to provide a witness Ferrers board in which the nineteen rows of Table~\ref{tab:singleton-classes} have different counts of avoiding placements. Such a suitable witness is the board with six rows of lengths $(6,6,6,6,6,4)$. The corresponding avoidance counts are listed in Table~\ref{tab:singleton-separation-six} in the Appendix.
\end{proof}

It turns out that the colored shape-Wilf classification of Table~\ref{tab:singleton-classes}, together with the known 
Wilf-classification~\cite{DMR2008} of signed patterns, provides all the needed information for colored Wilf-classification of these patterns. Recall that Mansour~\cite{mansour} has already shown that colored patterns of size 2 form a single colored Wilf class. 

Let us say that two colored patterns $\pi$ and $\sigma$ are \emph{symmetric} if $\pi$ can be turned into $\sigma$ by 
a series of horizontal, vertical and diagonal reflections and color-isomorphism. Clearly symmetric patterns are in the same colored Wilf class, though not necessarily in the same colored shape-Wilf class.

\begin{corollary}
Colored patterns of size 3 form two colored Wilf classes: one of them contains all the patterns symmetric to a pattern from Class VII in Table~\ref{tab:singleton-classes}, and the other contains the rest.
\end{corollary}

\begin{proof}
By inspection, any colored pattern of size three is either symmetric to a pattern from Class VII, or to a pattern from Class III in Table~\ref{tab:singleton-classes}. Thus, there can be no more than two colored Wilf classes. On the other hand, Dukes et al.~\cite{DMR2008}, working in the terminology of signed patterns, have shown that 2-colored patterns of size 3 already form two distinct 2-color Wilf classes. This implies that the two potential colored Wilf classes identified above are indeed distinct. 
\end{proof}

\section{Pairs of patterns of size 2}\label{sec:sets-of-patterns}

Let us now apply our methods to the classification of pairs of colored patterns of size~2. Again, we only need to consider one representative from each color-isomorphism class. Note that the color isomorphism must be applied to the whole pair of patterns, so we cannot limit ourselves to sets which only contain canonical patterns. To choose a representative for each color-isomorphism class, we shall use the following rules:
\begin{itemize}
    \item If a pair of patterns includes a coloring of 12 and a coloring of 21, we put the 12 first.
    \item If a pair of patterns includes two colorings of the same permutation and one of them is monochromatic, we put the monochromatic pattern first.
    \item We color the four entries by the same rule as for the canonical coloring of a permutation of size 4.
    \item We choose \patset{12/RB,12/BG} as representative over the color-isomorphic pair \patset{12/RB,12/GR}, and also \patset{21/RB,21/BG} over \patset{21/RB,21/GR}.
\end{itemize}
We will call a pair of colored patterns of size 2 that are colored by the above rules \emph{canonical}, and will henceforth use them as representatives of their color-isomorphism class. 

We begin by listing the identified colored shape-Wilf classes in Table~\ref{tab:two-color-pair-classes}. 

\begin{table}[ht]
\centering
\renewcommand{\arraystretch}{1.15}
\begin{tabular}[c]{cl}
\toprule
Class & Canonical pairs \\
\midrule
I &\patset{12/RR,12/RB}, \patset{12/RR,21/RR}, \patset{12/RB,12/RG}\\ \hline
II &
\begin{tabular}[t]{@{}l@{}}
\patset{12/RR,12/BR}, \patset{21/RR,21/RB}, \patset{21/RR,21/BR},\\ 
\patset{12/RB,21/RB}, \patset{12/RB,21/BR}, \patset{12/RB,12/GB}
\end{tabular}
\\\hline
III  &
\begin{tabular}[t]{@{}l@{}}
 \patset{12/RR,12/BB}, \patset{21/RR,21/BB}, \patset{12/RR,21/BB},\\
 \patset{12/RR,12/BG}, \patset{12/RB,21/GG}, \patset{12/RB,12/GY}
\end{tabular}
\\\hline
IV &  \patset{12/RB,12/BR}  \\\hline
V &  \patset{21/RB,21/BR}  \\\hline
VI  &
 \patset{12/RR,21/BR}\\\hline
VII & \patset{12/RB,21/RR}\\\hline
VIII & \patset{12/RR,21/RB}\\\hline
IX & \patset{12/RB,21/BB}\\\hline
X & \patset{12/RB,12/BG} \\\hline
XI & \patset{12/RB,21/GR} \\\hline
XII & \patset{12/RB,21/GB} \\\hline
XIII & \patset{12/RR,21/BG}, \patset{21/RR,21/BG}, \patset{12/RB,21/GY} \\\hline
XIV & \patset{12/RB,21/BG} \\\hline
XV & \patset{12/RB,21/RG} \\\hline
XVI & \patset{21/RB,21/BG} \\\hline
XVII & \patset{21/RB,21/GB} \\\hline
XVIII & \patset{21/RB,21/RG} \\\hline
XIX & \patset{21/RB,21/GY} \\
\bottomrule
\end{tabular}
\caption{The colored shape-Wilf equivalence classes of pairs of patterns of size~2, with a canonical representative of each color-isomorphism class.}
\label{tab:two-color-pair-classes}
\end{table}

Unlike the case of a single pattern of size 3, we cannot explain all the equivalences in Table~\ref{tab:two-color-pair-classes} by the general properties established in Section~\ref{sec-general}. Instead, we have to deal with two cases individually in the next two lemmas.

\begin{lemma}\label{lem-spec1}
    The two pairs $\patset{12/RR,12/RB}$ and \patset{12/RR,21/RR} are colored shape-Wilf equivalent. 
\end{lemma}

\begin{proof}
    Let $P= \patset{12/RR,12/RB}$ and $Q= \patset{12/RR,21/RR}$. We will show that $P$ and $Q$ are 2-colored shape-Wilf equivalent.

    Let $F$ be a Ferrers board of order~$n$. Let us say that a cell $c\in F$ is \emph{admissible} if $F$ has a placement in which $c$ is a 1-cell. For an admissible cell $c$, let $F/c$ be the subboard of $F$ obtained by
    deleting the row and column containing $c$ and closing the gaps. 
    
    By induction on $n$, we will show that $F$ has the same number of $P$-avoiding and $Q$-avoiding 2-colored placements, i.e.,
\begin{equation}
|S_F^2(P)|=|S_F^2(Q)|.\label{eq-ident}
\end{equation}
    
    This is trivially true for $n=1$, so let us assume that $n\ge 2$ and proceed with the induction step. 
    Let $R$ be the set of row-indices $j$ such that the cell $(1,j)$ (i.e., the cell in the first column
    and $j$-th row) is admissible. If $R$ is empty, then $F$ has no placement and \eqref{eq-ident} holds vacuously.
    Suppose from now on that $R\neq\emptyset$ and notice that if $R$ contains a row index $j<n$, then it also contains all the higher indices $j+1,j+2,\dotsc,n$.
    
    For each $j \in R$ we let $F^{(j)}$ denote the Ferrers board $F/(1,j)$. 
    Let us look at a fixed $P$-avoiding placement on some $F^{(k)}$, and see how we can extend it to a $P$-avoiding placement on~$F$. We first reinsert the first column and the $k$-th row, and then we can add a blue rook at the $(1,k)$ position, which would yield a $P$-avoiding placement in~$F$. Adding a red rook, however, is only possible when $n=k$; otherwise, we would get either $\pat{12/RR}$ or $\pat{12/RB}$ in the placement, so in this case we actually add back the first column and the top row. 

    So, we have the following recurrence relation:

\[     |S_F^2(P)| = \sum_{k\in R} |S_{F^{(k)}}^2(P)|+ |S_{F^{(n)}}^2(P)|\]

    Now, let us look at a fixed $Q$-avoiding placement in some $F^{(k)}$, and extend it into a placement on~$F$ by adding a rook to the position $(1,k)$. If we are adding a blue rook, no occurrence of $Q$ will be created, which contributes $\sum_{k\in R} |S_{F^{(k)}}^2(Q)|$ placements with a blue rook in the first column. 
    
To account for $Q$-avoiding placements on $F$ with a red rook in the first column, we will show that they are 
in bijection with $Q$-avoiding placements in~$F^{(n)}$. Let $d$ be the length of the top row of~$F$. Fix a $Q$-avoiding placement $T$ in~$F^{(n)}$.

If in the first $d-1$ columns of $T$, all the rooks are blue, then we can add a red rook at position 
$(1,n)$ to get a $Q$-avoiding placement in~$F$. If not, then we have exactly one red rook in those $d-1$ columns; otherwise we would either get $\pat{12/RR}$ or $\pat{21/RR}$ in~$T$. Suppose the red rook is in row $j$ and column $i$ with respect to the original board $F$ (so in column $i-1$ in $F^{(n)}$). To get a placement in $F$, move the red rook in row $j$ to the first column, and in the topmost row, add a blue rook in column~$i$.

Since the columns $1,\dots,d$ in $F$ have the same height, transferring the red rook between columns $i$ and $1$ in 
the same row preserves its compatibility with every other rook. In particular, it does not create any occurrence of 
$\pat{12/RR}$ or $\pat{21/RR}$, and the resulting placement in $F$ avoids~$Q$.
    With a similar argument, we have the inverse of this procedure; thus, we have a bijection between the $Q$-avoiding 2-colored placements of $F^{(n)}$ and the $Q$-avoiding 2-colored placements in $F$ that have a red rook in the first column.
    
    Therefore, we have the recurrence relation:
\[ 
|S_F^2(Q)| = \sum_{k\in R} |S_{F^{(k)}}^2(Q)|+ |S_{F^{(n)}}^2(Q)|.
\]
     Therefore, by induction, we have $|S_F^2(P)| = |S_F^2(Q)|$.
\end{proof}

We proceed with our second single-purpose lemma.

\begin{lemma}\label{lem:C2-mixed}
The set $P=\patset{12/RR,12/BR}$ is colored shape-Wilf equivalent to the set $Q=\patset{12/RB,21/RB}$.
\end{lemma}
\begin{proof}
We apply a similar inductive strategy as in the proof of Lemma~\ref{lem-spec1}, but this time, we will eliminate the top row and one of its intersecting columns. 

Let $\lambda=(\lambda_1\ge\dots\ge\lambda_n)$
be the sequence of row lengths of a Ferrers board~$F$, and put
$d=\lambda_n$. We may assume that $F$ has at least two rows and two columns, since the patterns we are considering 
have order two. We may also assume that $F$ admits a placement.

Let $F^{-}$ be the Ferrers board with row lengths
$(\lambda_1-1,\dots,\lambda_{n-1}-1).$
Every 2-colored placement in $F$ is obtained uniquely from a placement in
$F^{-}$ by inserting a new column in one of the $d$ possible positions and
placing a rook in this column and the new top row.

For every $P$-avoiding placement in $F^{-}$, a new blue rook
may be inserted in any of the $d$ positions of the top row.
A new red rook must have no compatible rook South-West
of it, so it may only be inserted in the leftmost position.
Thus there are exactly $d+1$ valid insertions, and 
\[
|S_F^2(P)|=(d+1)|S_{F^{-}}^2(P)|.
\]

Now consider $Q=\patset{12/RB,21/RB}$ and a $Q$-avoiding placement $T$ in~$F^{-}$.
The leftmost $d-1$ columns have full height, so their rooks are pairwise compatible.
Since $T$ avoids $Q$, all the red rooks in these columns are to the right of any blue rooks there.

Suppose there are $q$ blue rooks overall in the leftmost $d-1$ columns of~$T$. Necessarily, 
they occupy columns $1,2,\dotsc,q$ in~$F^{-}$. 

Now insert into $T$ a new top row of length $d$ and a new column intersecting it, and place a rook at their intersection, to turn $T$ into a placement in~$F$.
The new rook is compatible precisely with the rooks
from the leftmost $d-1$ columns of~$T$. To avoid $Q$, a new blue rook has $q+1$ possible positions, while
a new red rook has $d-q$ possible positions.
Therefore, there are again $(q+1)+(d-q)=d+1$ valid insertions, and we get 
\[
|S_F^2(Q)|=(d+1)|S_{F^{-}}^2(Q)|.
\]

The induction hypothesis therefore gives
\[
  \patset{12/RR,12/BR}\swilf[2]\patset{12/RB,21/RB}.\qedhere
\]
\end{proof}

Now we are ready to complete the classification of pairs of patterns of size~2.

\begin{theorem}
There are nineteen colored shape-Wilf equivalence classes among the pairs of colored patterns of size 2. They are listed in Table~\ref{tab:two-color-pair-classes}.
\end{theorem}
\begin{proof}
We begin by showing that the claimed classes indeed contain equivalent pairs. 

\emph{Class I.} The equivalence of \patset{12/RR,12/RB} and \patset{12/RR,21/RR} is shown in Lemma~\ref{lem-spec1}. For the pair \patset{12/RB,12/RG}, note that it can be written as $1_R\oplus\{1_B,1_G\}$, which by Lemma~\ref{lem:tail-recoloring} is equivalent to $1_R\oplus\{1_B,1_R\}$, which is the pair \patset{12/RR,12/RB}.

\emph{Class II.} 
We first deal with $P_1=\patset{12/RR,12/BR}$, $P_2=\patset{21/RR,21/RB}$ and 
$P_3=\patset{21/RR,21/BR}$. Observe that a 2-colored placement avoids $P_1$ if and only if every red rook is 
SW-extreme. Likewise, it avoids $P_2$ if and only if every red rook is SE-extreme, and it avoids $P_3$ if and only if 
every red rook is NW-extreme. A given uncolored placement $T$ therefore admits exactly $2^{\SW(T)}$ 2-colorings 
avoiding $P_1$, since we have two choices of colors for each SW-extreme rook, and the remaining rooks are blue. Summing this over all uncolored placements $T$ of a given Ferrers board $F$ gives us the number of $P_1$-avoiding 2-colored placements in~$F$. An analogous accounting works for $P_2$ and $P_3$, with the statistics $\SE$ and $\NW$ in place of~$\SW$. By Theorem~\ref{thm-equi}, the three counts of avoiding placements are the same for every Ferrers board~$F$.

The equivalence of $P_1=\patset{12/RR,12/BR}$ and $P_4=\patset{12/RB,21/RB}$ is shown in Lemma~\ref{lem:C2-mixed}. Using Observation~\ref{obs-transpose} and noting that $P_1^t=P_1$, we see that $P_1$ is also equivalent to 
$P_4^t=\patset{12/RB,21/BR}$.

To show that $P_5=\patset{12/RB,12/GB}$ is equivalent to $P_1$, observe that $P_5=\{1_R,1_G\}\oplus 1_B$ and $P_1=\{1_R,1_B\}\oplus 1_R$. Their equivalence then follows from Corollary~\ref{cor:two-changes}.

\emph{Class III.} Note that this class contains precisely the pairs in which the two patterns do not have any color in common and each pattern is either a monochromatic coloring of 21 or an arbitrary coloring of 12. From the uncolored shape-Wilf equivalence of 12 and 21, we obtain that any monochromatic coloring of 12 is colored shape-Wilf equivalent to any monochromatic coloring of 21. A repeated application of Lemma~\ref{lem:pattern-addition} then yields the colored shape-Wilf equivalence of the pairs \patset{12/RR,12/BB}, \patset{12/RR,21/BB} and \patset{21/RR,21/BB}. 

By Proposition~\ref{prop:1n}, we know that all colorings of 12 are colored shape-Wilf equivalent. Using this fact in conjunction with Lemma~\ref{lem:pattern-addition} and with color-isomorphism if needed, we obtain that the remaining three pairs also belong to the same equivalence class as the first three. 

\emph{Class XIII.} This class again contains pairs of patterns on disjoint color sets, but this time the second 
pattern in each pair is color-isomorphic to \pat{21/RB}. After applying a color-isomorphism that makes the second patterns in the three pairs identical, their equivalence follows directly from Lemma~\ref{lem:pattern-addition}.

To argue that no other pairs are equivalent, it is again enough to enumerate their avoiders in a single Ferrers board of order 6; see Table~\ref{tab:pair-class-separation} in the Appendix.
\end{proof}

What does this tell us about colored Wilf-equivalence for these pairs of patterns? If we combine those pairs that can be related by symmetry, the nineteen classes of Table~\ref{tab:two-color-pair-classes} further merge into the following six superclasses:
\begin{enumerate}
    \item The union of Classes I, II, XVII and XVIII.
    \item The union of Classes III, XIII and XIX (note that this contains exactly the pairs where the two patterns have disjoint color sets).
    \item The union of Classes IV and V.
    \item The union of Classes VI, VII, VIII and IX.
    \item The union of Classes X and XVI.
    \item The union of Classes XI, XII, XIV and XV.
\end{enumerate}
This recovers exactly the six colored Wilf-equivalence classes that have been already determined by Mansour~\cite[Table 1]{mansour}. Note that only the first two superclasses above include more than one orbit under symmetry and color-isomorphism.

\section{Conclusion} 
Having succeeded in classifying colored patterns of size up to 3 and pairs of patterns of size 2, an obvious follow-up is to expand this to larger patterns or larger sets, in the hope that perhaps further general relationships will emerge. So far, we have made no effort to move in this direction, or even to look for plausible conjectures to make.

\section*{Acknowledgments}
This work was supported by the Student Faculty Grant (SFG)
of the Faculty of Mathematics and Physics,
Charles University, Prague, Czech Republic.

\section*{AI Disclosure} 
The bijective argument in the proof of Lemma~\ref{lem:C2-mixed} was proposed by Open\-AI's ChatGPT 6.0 Sol model; the two human authors have verified the correctness and originality of the argument. GPT 5.6 Sol in high mode was used to generate Figure~\ref{fig:q-coloring}, to typeset Tables~\ref{tab:two-color-pair-classes}, \ref{tab:singleton-separation-six} and \ref{tab:pair-class-separation} from equivalences/avoidance counts supplied by the authors in a text file, and to find typos or mistakes in various draft versions of the paper. The authors have also used GPT 5.6 Sol in medium mode and the web version of ChatGPT to assist in searching for relevant references. 

The two authors reviewed all AI-assisted contributions, and they have written the whole text of this paper 
by themselves. They take full responsibility for the correctness, originality, attribution,
and presentation of the paper, including its references.

\bibliographystyle{plain}
\bibliography{references}

\section{Appendix}

\begin{table}[H]
\centering
\small
\begin{tabular}{clrr}
\toprule
Class & Representative $p$ & $|S_F^2(p)|$ & $|S_F^3(p)|$ \\
\midrule
Length 1 & $\pat{1/R}$ & 480 & 30720 \\
I & $\pat{12/RR}$ & 9875 & 205829 \\
II & $\pat{21/RB}$ & 9085 & 199789 \\
III & $\pat{123/RRR}$ & 25524 & 329162 \\
IV & $\pat{321/RRB}$ & 25211 & 328239 \\
V & $\pat{321/RBR}$ & 25120 & 328078 \\
VI & $\pat{321/RBB}$ & 25147 & 328115 \\
VII & $\pat{213/RBR}$ & 25480 & 329084 \\
VIII & $\pat{231/RRR}$ & 25226 & 328264 \\
IX & $\pat{132/RBB}$ & 25538 & 329186 \\
X & $\pat{132/RRB}$ & 25218 & 328252 \\
XI & $\pat{231/RRB}$ & 25175 & 328173 \\
XII & $\pat{231/RBR}$ & 25190 & 328204 \\
XIII & $\pat{231/RBB}$ & 25199 & 328213 \\
XIV & $\pat{312/RRB}$ & 25144 & 328118 \\
XV & $\pat{312/RBR}$ & 25220 & 328254 \\
XVI & $\pat{312/RBB}$ & 25131 & 328089 \\
XVII & $\pat{321/RBG}$ & --- & 328012 \\
XVIII & $\pat{231/RBG}$ & --- & 328138 \\
XIX & $\pat{312/RBG}$ & --- & 328052 \\
\bottomrule
\end{tabular}
\caption{Avoidance counts for representatives of equivalence classes of singleton patterns of length at most three, using at most three colors on board $F=(6,6,6,6,6,4)$.}
\label{tab:singleton-separation-six}
\end{table}

\begin{table}[H]
\centering
\small
\renewcommand{\arraystretch}{1.08}
\begin{tabular}{clr}
\toprule
Class & Representative $P$ & $|S_F^4(P)|$ \\
\midrule
I & $\patset{12/RR,12/RB}$ & 563760 \\
II & $\patset{12/RR,12/BR}$ & 538704 \\
III & $\patset{12/RR,12/BB}$ & 511192 \\
IV & $\patset{12/RB,12/BR}$ & 524750 \\
V & $\patset{21/RB,21/BR}$ & 524842 \\
VI & $\patset{12/RR,21/BR}$ & 532535 \\
VII & $\patset{12/RB,21/RR}$ & 558429 \\
VIII & $\patset{12/RR,21/RB}$ & 532041 \\
IX & $\patset{12/RB,21/BB}$ & 531408 \\
X & $\patset{12/RB,12/BG}$ & 517970 \\
XI & $\patset{12/RB,21/GR}$ & 526508 \\
XII & $\patset{12/RB,21/GB}$ & 511404 \\
XIII & $\patset{12/RR,21/BG}$ & 495846 \\
XIV & $\patset{12/RB,21/BG}$ & 510885 \\
XV & $\patset{12/RB,21/RG}$ & 525520 \\
XVI & $\patset{21/RB,21/BG}$ & 504422 \\
XVII & $\patset{21/RB,21/GB}$ & 515856 \\
XVIII & $\patset{21/RB,21/RG}$ & 515984 \\
XIX & $\patset{21/RB,21/GY}$ & 480572 \\
\bottomrule
\end{tabular}
\caption{Avoidance counts for representatives
of the 19  classes of two-element sets of patterns
of length two.
The Ferrers board is $F=(6,6,6,6,4,3)$,
and $R,B,G,Y$ denote distinct colors.}
\label{tab:pair-class-separation}
\end{table}

\begingroup
\small
\renewcommand{\arraystretch}{1.08}
\setlength{\tabcolsep}{4pt}
\setlength{\LTcapwidth}{\linewidth}
\setlength{\LTleft}{0pt}
\setlength{\LTright}{0pt}
\begin{longtable}{@{}
  >{\centering\arraybackslash}p{3em}
  >{\raggedright\arraybackslash}p{\dimexpr\linewidth-3em-2\tabcolsep\relax}
  @{}}
\caption{All the 19 classes of two-element bases
of length-two patterns over four colors.}
\label{tab:all-four-color-pair-bases}\\
\toprule
Class & Members \\
\midrule
\endfirsthead
\multicolumn{2}{c}{\tablename~\thetable\ (continued)}\\
\toprule
Class & Members \\
\midrule
\endhead
\midrule
\multicolumn{2}{r}{\footnotesize Continued on the next page}\\
\endfoot
\bottomrule
\endlastfoot
$\mathcal C_{1}$ & $\patset{12/RR,12/RB}$,\allowbreak\ $\patset{12/RR,12/RG}$,\allowbreak\ $\patset{12/RR,12/RY}$,\allowbreak\ $\patset{12/RR,21/RR}$,\allowbreak\ $\patset{12/RB,12/RG}$,\allowbreak\ $\patset{12/RB,12/RY}$ \\
$\mathcal C_{1}$ & $\patset{12/RG,12/RY}$,\allowbreak\ $\patset{12/BR,12/BB}$,\allowbreak\ $\patset{12/BR,12/BG}$,\allowbreak\ $\patset{12/BR,12/BY}$,\allowbreak\ $\patset{12/BB,12/BG}$,\allowbreak\ $\patset{12/BB,12/BY}$ \\
$\mathcal C_{1}$ & $\patset{12/BB,21/BB}$,\allowbreak\ $\patset{12/BG,12/BY}$,\allowbreak\ $\patset{12/GR,12/GB}$,\allowbreak\ $\patset{12/GR,12/GG}$,\allowbreak\ $\patset{12/GR,12/GY}$,\allowbreak\ $\patset{12/GB,12/GG}$ \\
$\mathcal C_{1}$ & $\patset{12/GB,12/GY}$,\allowbreak\ $\patset{12/GG,12/GY}$,\allowbreak\ $\patset{12/GG,21/GG}$,\allowbreak\ $\patset{12/YR,12/YB}$,\allowbreak\ $\patset{12/YR,12/YG}$,\allowbreak\ $\patset{12/YR,12/YY}$ \\
$\mathcal C_{1}$ & $\patset{12/YB,12/YG}$,\allowbreak\ $\patset{12/YB,12/YY}$,\allowbreak\ $\patset{12/YG,12/YY}$,\allowbreak\ $\patset{12/YY,21/YY}$ \\
\midrule
$\mathcal C_{2}$ & $\patset{12/RR,12/BR}$,\allowbreak\ $\patset{12/RR,12/GR}$,\allowbreak\ $\patset{12/RR,12/YR}$,\allowbreak\ $\patset{12/RB,12/BB}$,\allowbreak\ $\patset{12/RB,12/GB}$,\allowbreak\ $\patset{12/RB,12/YB}$ \\
$\mathcal C_{2}$ & $\patset{12/RB,21/RB}$,\allowbreak\ $\patset{12/RB,21/BR}$,\allowbreak\ $\patset{12/RG,12/BG}$,\allowbreak\ $\patset{12/RG,12/GG}$,\allowbreak\ $\patset{12/RG,12/YG}$,\allowbreak\ $\patset{12/RG,21/RG}$ \\
$\mathcal C_{2}$ & $\patset{12/RG,21/GR}$,\allowbreak\ $\patset{12/RY,12/BY}$,\allowbreak\ $\patset{12/RY,12/GY}$,\allowbreak\ $\patset{12/RY,12/YY}$,\allowbreak\ $\patset{12/RY,21/RY}$,\allowbreak\ $\patset{12/RY,21/YR}$ \\
$\mathcal C_{2}$ & $\patset{12/BR,12/GR}$,\allowbreak\ $\patset{12/BR,12/YR}$,\allowbreak\ $\patset{12/BR,21/RB}$,\allowbreak\ $\patset{12/BR,21/BR}$,\allowbreak\ $\patset{12/BB,12/GB}$,\allowbreak\ $\patset{12/BB,12/YB}$ \\
$\mathcal C_{2}$ & $\patset{12/BG,12/GG}$,\allowbreak\ $\patset{12/BG,12/YG}$,\allowbreak\ $\patset{12/BG,21/BG}$,\allowbreak\ $\patset{12/BG,21/GB}$,\allowbreak\ $\patset{12/BY,12/GY}$,\allowbreak\ $\patset{12/BY,12/YY}$ \\
$\mathcal C_{2}$ & $\patset{12/BY,21/BY}$,\allowbreak\ $\patset{12/BY,21/YB}$,\allowbreak\ $\patset{12/GR,12/YR}$,\allowbreak\ $\patset{12/GR,21/RG}$,\allowbreak\ $\patset{12/GR,21/GR}$,\allowbreak\ $\patset{12/GB,12/YB}$ \\
$\mathcal C_{2}$ & $\patset{12/GB,21/BG}$,\allowbreak\ $\patset{12/GB,21/GB}$,\allowbreak\ $\patset{12/GG,12/YG}$,\allowbreak\ $\patset{12/GY,12/YY}$,\allowbreak\ $\patset{12/GY,21/GY}$,\allowbreak\ $\patset{12/GY,21/YG}$ \\
$\mathcal C_{2}$ & $\patset{12/YR,21/RY}$,\allowbreak\ $\patset{12/YR,21/YR}$,\allowbreak\ $\patset{12/YB,21/BY}$,\allowbreak\ $\patset{12/YB,21/YB}$,\allowbreak\ $\patset{12/YG,21/GY}$,\allowbreak\ $\patset{12/YG,21/YG}$ \\
$\mathcal C_{2}$ & $\patset{21/RR,21/RB}$,\allowbreak\ $\patset{21/RR,21/RG}$,\allowbreak\ $\patset{21/RR,21/RY}$,\allowbreak\ $\patset{21/RR,21/BR}$,\allowbreak\ $\patset{21/RR,21/GR}$,\allowbreak\ $\patset{21/RR,21/YR}$ \\
$\mathcal C_{2}$ & $\patset{21/RB,21/BB}$,\allowbreak\ $\patset{21/RG,21/GG}$,\allowbreak\ $\patset{21/RY,21/YY}$,\allowbreak\ $\patset{21/BR,21/BB}$,\allowbreak\ $\patset{21/BB,21/BG}$,\allowbreak\ $\patset{21/BB,21/BY}$ \\
$\mathcal C_{2}$ & $\patset{21/BB,21/GB}$,\allowbreak\ $\patset{21/BB,21/YB}$,\allowbreak\ $\patset{21/BG,21/GG}$,\allowbreak\ $\patset{21/BY,21/YY}$,\allowbreak\ $\patset{21/GR,21/GG}$,\allowbreak\ $\patset{21/GB,21/GG}$ \\
$\mathcal C_{2}$ & $\patset{21/GG,21/GY}$,\allowbreak\ $\patset{21/GG,21/YG}$,\allowbreak\ $\patset{21/GY,21/YY}$,\allowbreak\ $\patset{21/YR,21/YY}$,\allowbreak\ $\patset{21/YB,21/YY}$,\allowbreak\ $\patset{21/YG,21/YY}$ \\
\midrule
$\mathcal C_{3}$ & $\patset{12/RR,12/BB}$,\allowbreak\ $\patset{12/RR,12/BG}$,\allowbreak\ $\patset{12/RR,12/BY}$,\allowbreak\ $\patset{12/RR,12/GB}$,\allowbreak\ $\patset{12/RR,12/GG}$,\allowbreak\ $\patset{12/RR,12/GY}$ \\
$\mathcal C_{3}$ & $\patset{12/RR,12/YB}$,\allowbreak\ $\patset{12/RR,12/YG}$,\allowbreak\ $\patset{12/RR,12/YY}$,\allowbreak\ $\patset{12/RR,21/BB}$,\allowbreak\ $\patset{12/RR,21/GG}$,\allowbreak\ $\patset{12/RR,21/YY}$ \\
$\mathcal C_{3}$ & $\patset{12/RB,12/GG}$,\allowbreak\ $\patset{12/RB,12/GY}$,\allowbreak\ $\patset{12/RB,12/YG}$,\allowbreak\ $\patset{12/RB,12/YY}$,\allowbreak\ $\patset{12/RB,21/GG}$,\allowbreak\ $\patset{12/RB,21/YY}$ \\
$\mathcal C_{3}$ & $\patset{12/RG,12/BB}$,\allowbreak\ $\patset{12/RG,12/BY}$,\allowbreak\ $\patset{12/RG,12/YB}$,\allowbreak\ $\patset{12/RG,12/YY}$,\allowbreak\ $\patset{12/RG,21/BB}$,\allowbreak\ $\patset{12/RG,21/YY}$ \\
$\mathcal C_{3}$ & $\patset{12/RY,12/BB}$,\allowbreak\ $\patset{12/RY,12/BG}$,\allowbreak\ $\patset{12/RY,12/GB}$,\allowbreak\ $\patset{12/RY,12/GG}$,\allowbreak\ $\patset{12/RY,21/BB}$,\allowbreak\ $\patset{12/RY,21/GG}$ \\
$\mathcal C_{3}$ & $\patset{12/BR,12/GG}$,\allowbreak\ $\patset{12/BR,12/GY}$,\allowbreak\ $\patset{12/BR,12/YG}$,\allowbreak\ $\patset{12/BR,12/YY}$,\allowbreak\ $\patset{12/BR,21/GG}$,\allowbreak\ $\patset{12/BR,21/YY}$ \\
$\mathcal C_{3}$ & $\patset{12/BB,12/GR}$,\allowbreak\ $\patset{12/BB,12/GG}$,\allowbreak\ $\patset{12/BB,12/GY}$,\allowbreak\ $\patset{12/BB,12/YR}$,\allowbreak\ $\patset{12/BB,12/YG}$,\allowbreak\ $\patset{12/BB,12/YY}$ \\
$\mathcal C_{3}$ & $\patset{12/BB,21/RR}$,\allowbreak\ $\patset{12/BB,21/GG}$,\allowbreak\ $\patset{12/BB,21/YY}$,\allowbreak\ $\patset{12/BG,12/YR}$,\allowbreak\ $\patset{12/BG,12/YY}$,\allowbreak\ $\patset{12/BG,21/RR}$ \\
$\mathcal C_{3}$ & $\patset{12/BG,21/YY}$,\allowbreak\ $\patset{12/BY,12/GR}$,\allowbreak\ $\patset{12/BY,12/GG}$,\allowbreak\ $\patset{12/BY,21/RR}$,\allowbreak\ $\patset{12/BY,21/GG}$,\allowbreak\ $\patset{12/GR,12/YB}$ \\
$\mathcal C_{3}$ & $\patset{12/GR,12/YY}$,\allowbreak\ $\patset{12/GR,21/BB}$,\allowbreak\ $\patset{12/GR,21/YY}$,\allowbreak\ $\patset{12/GB,12/YR}$,\allowbreak\ $\patset{12/GB,12/YY}$,\allowbreak\ $\patset{12/GB,21/RR}$ \\
$\mathcal C_{3}$ & $\patset{12/GB,21/YY}$,\allowbreak\ $\patset{12/GG,12/YR}$,\allowbreak\ $\patset{12/GG,12/YB}$,\allowbreak\ $\patset{12/GG,12/YY}$,\allowbreak\ $\patset{12/GG,21/RR}$,\allowbreak\ $\patset{12/GG,21/BB}$ \\
$\mathcal C_{3}$ & $\patset{12/GG,21/YY}$,\allowbreak\ $\patset{12/GY,21/RR}$,\allowbreak\ $\patset{12/GY,21/BB}$,\allowbreak\ $\patset{12/YR,21/BB}$,\allowbreak\ $\patset{12/YR,21/GG}$,\allowbreak\ $\patset{12/YB,21/RR}$ \\
$\mathcal C_{3}$ & $\patset{12/YB,21/GG}$,\allowbreak\ $\patset{12/YG,21/RR}$,\allowbreak\ $\patset{12/YG,21/BB}$,\allowbreak\ $\patset{12/YY,21/RR}$,\allowbreak\ $\patset{12/YY,21/BB}$,\allowbreak\ $\patset{12/YY,21/GG}$ \\
$\mathcal C_{3}$ & $\patset{21/RR,21/BB}$,\allowbreak\ $\patset{21/RR,21/GG}$,\allowbreak\ $\patset{21/RR,21/YY}$,\allowbreak\ $\patset{21/BB,21/GG}$,\allowbreak\ $\patset{21/BB,21/YY}$,\allowbreak\ $\patset{21/GG,21/YY}$ \\
\midrule
$\mathcal C_{4}$ & $\patset{12/RB,12/BR}$,\allowbreak\ $\patset{12/RG,12/GR}$,\allowbreak\ $\patset{12/RY,12/YR}$,\allowbreak\ $\patset{12/BG,12/GB}$,\allowbreak\ $\patset{12/BY,12/YB}$,\allowbreak\ $\patset{12/GY,12/YG}$ \\
\midrule
$\mathcal C_{5}$ & $\patset{21/RB,21/BR}$,\allowbreak\ $\patset{21/RG,21/GR}$,\allowbreak\ $\patset{21/RY,21/YR}$,\allowbreak\ $\patset{21/BG,21/GB}$,\allowbreak\ $\patset{21/BY,21/YB}$,\allowbreak\ $\patset{21/GY,21/YG}$ \\
\midrule
$\mathcal C_{6}$ & $\patset{12/RR,21/BR}$,\allowbreak\ $\patset{12/RR,21/GR}$,\allowbreak\ $\patset{12/RR,21/YR}$,\allowbreak\ $\patset{12/BB,21/RB}$,\allowbreak\ $\patset{12/BB,21/GB}$,\allowbreak\ $\patset{12/BB,21/YB}$ \\
$\mathcal C_{6}$ & $\patset{12/GG,21/RG}$,\allowbreak\ $\patset{12/GG,21/BG}$,\allowbreak\ $\patset{12/GG,21/YG}$,\allowbreak\ $\patset{12/YY,21/RY}$,\allowbreak\ $\patset{12/YY,21/BY}$,\allowbreak\ $\patset{12/YY,21/GY}$ \\
\midrule
$\mathcal C_{7}$ & $\patset{12/RB,21/RR}$,\allowbreak\ $\patset{12/RG,21/RR}$,\allowbreak\ $\patset{12/RY,21/RR}$,\allowbreak\ $\patset{12/BR,21/BB}$,\allowbreak\ $\patset{12/BG,21/BB}$,\allowbreak\ $\patset{12/BY,21/BB}$ \\
$\mathcal C_{7}$ & $\patset{12/GR,21/GG}$,\allowbreak\ $\patset{12/GB,21/GG}$,\allowbreak\ $\patset{12/GY,21/GG}$,\allowbreak\ $\patset{12/YR,21/YY}$,\allowbreak\ $\patset{12/YB,21/YY}$,\allowbreak\ $\patset{12/YG,21/YY}$ \\
\midrule
$\mathcal C_{8}$ & $\patset{12/RR,21/RB}$,\allowbreak\ $\patset{12/RR,21/RG}$,\allowbreak\ $\patset{12/RR,21/RY}$,\allowbreak\ $\patset{12/BB,21/BR}$,\allowbreak\ $\patset{12/BB,21/BG}$,\allowbreak\ $\patset{12/BB,21/BY}$ \\
$\mathcal C_{8}$ & $\patset{12/GG,21/GR}$,\allowbreak\ $\patset{12/GG,21/GB}$,\allowbreak\ $\patset{12/GG,21/GY}$,\allowbreak\ $\patset{12/YY,21/YR}$,\allowbreak\ $\patset{12/YY,21/YB}$,\allowbreak\ $\patset{12/YY,21/YG}$ \\
\midrule
$\mathcal C_{9}$ & $\patset{12/RB,21/BB}$,\allowbreak\ $\patset{12/RG,21/GG}$,\allowbreak\ $\patset{12/RY,21/YY}$,\allowbreak\ $\patset{12/BR,21/RR}$,\allowbreak\ $\patset{12/BG,21/GG}$,\allowbreak\ $\patset{12/BY,21/YY}$ \\
$\mathcal C_{9}$ & $\patset{12/GR,21/RR}$,\allowbreak\ $\patset{12/GB,21/BB}$,\allowbreak\ $\patset{12/GY,21/YY}$,\allowbreak\ $\patset{12/YR,21/RR}$,\allowbreak\ $\patset{12/YB,21/BB}$,\allowbreak\ $\patset{12/YG,21/GG}$ \\
\midrule
$\mathcal C_{10}$ & $\patset{12/RB,12/BG}$,\allowbreak\ $\patset{12/RB,12/BY}$,\allowbreak\ $\patset{12/RB,12/GR}$,\allowbreak\ $\patset{12/RB,12/YR}$,\allowbreak\ $\patset{12/RG,12/BR}$,\allowbreak\ $\patset{12/RG,12/GB}$ \\
$\mathcal C_{10}$ & $\patset{12/RG,12/GY}$,\allowbreak\ $\patset{12/RG,12/YR}$,\allowbreak\ $\patset{12/RY,12/BR}$,\allowbreak\ $\patset{12/RY,12/GR}$,\allowbreak\ $\patset{12/RY,12/YB}$,\allowbreak\ $\patset{12/RY,12/YG}$ \\
$\mathcal C_{10}$ & $\patset{12/BR,12/GB}$,\allowbreak\ $\patset{12/BR,12/YB}$,\allowbreak\ $\patset{12/BG,12/GR}$,\allowbreak\ $\patset{12/BG,12/GY}$,\allowbreak\ $\patset{12/BG,12/YB}$,\allowbreak\ $\patset{12/BY,12/GB}$ \\
$\mathcal C_{10}$ & $\patset{12/BY,12/YR}$,\allowbreak\ $\patset{12/BY,12/YG}$,\allowbreak\ $\patset{12/GR,12/YG}$,\allowbreak\ $\patset{12/GB,12/YG}$,\allowbreak\ $\patset{12/GY,12/YR}$,\allowbreak\ $\patset{12/GY,12/YB}$ \\
\midrule
$\mathcal C_{11}$ & $\patset{12/RB,21/GR}$,\allowbreak\ $\patset{12/RB,21/YR}$,\allowbreak\ $\patset{12/RG,21/BR}$,\allowbreak\ $\patset{12/RG,21/YR}$,\allowbreak\ $\patset{12/RY,21/BR}$,\allowbreak\ $\patset{12/RY,21/GR}$ \\
$\mathcal C_{11}$ & $\patset{12/BR,21/GB}$,\allowbreak\ $\patset{12/BR,21/YB}$,\allowbreak\ $\patset{12/BG,21/RB}$,\allowbreak\ $\patset{12/BG,21/YB}$,\allowbreak\ $\patset{12/BY,21/RB}$,\allowbreak\ $\patset{12/BY,21/GB}$ \\
$\mathcal C_{11}$ & $\patset{12/GR,21/BG}$,\allowbreak\ $\patset{12/GR,21/YG}$,\allowbreak\ $\patset{12/GB,21/RG}$,\allowbreak\ $\patset{12/GB,21/YG}$,\allowbreak\ $\patset{12/GY,21/RG}$,\allowbreak\ $\patset{12/GY,21/BG}$ \\
$\mathcal C_{11}$ & $\patset{12/YR,21/BY}$,\allowbreak\ $\patset{12/YR,21/GY}$,\allowbreak\ $\patset{12/YB,21/RY}$,\allowbreak\ $\patset{12/YB,21/GY}$,\allowbreak\ $\patset{12/YG,21/RY}$,\allowbreak\ $\patset{12/YG,21/BY}$ \\
\midrule
$\mathcal C_{12}$ & $\patset{12/RB,21/GB}$,\allowbreak\ $\patset{12/RB,21/YB}$,\allowbreak\ $\patset{12/RG,21/BG}$,\allowbreak\ $\patset{12/RG,21/YG}$,\allowbreak\ $\patset{12/RY,21/BY}$,\allowbreak\ $\patset{12/RY,21/GY}$ \\
$\mathcal C_{12}$ & $\patset{12/BR,21/GR}$,\allowbreak\ $\patset{12/BR,21/YR}$,\allowbreak\ $\patset{12/BG,21/RG}$,\allowbreak\ $\patset{12/BG,21/YG}$,\allowbreak\ $\patset{12/BY,21/RY}$,\allowbreak\ $\patset{12/BY,21/GY}$ \\
$\mathcal C_{12}$ & $\patset{12/GR,21/BR}$,\allowbreak\ $\patset{12/GR,21/YR}$,\allowbreak\ $\patset{12/GB,21/RB}$,\allowbreak\ $\patset{12/GB,21/YB}$,\allowbreak\ $\patset{12/GY,21/RY}$,\allowbreak\ $\patset{12/GY,21/BY}$ \\
$\mathcal C_{12}$ & $\patset{12/YR,21/BR}$,\allowbreak\ $\patset{12/YR,21/GR}$,\allowbreak\ $\patset{12/YB,21/RB}$,\allowbreak\ $\patset{12/YB,21/GB}$,\allowbreak\ $\patset{12/YG,21/RG}$,\allowbreak\ $\patset{12/YG,21/BG}$ \\
\midrule
$\mathcal C_{13}$ & $\patset{12/RR,21/BG}$,\allowbreak\ $\patset{12/RR,21/BY}$,\allowbreak\ $\patset{12/RR,21/GB}$,\allowbreak\ $\patset{12/RR,21/GY}$,\allowbreak\ $\patset{12/RR,21/YB}$,\allowbreak\ $\patset{12/RR,21/YG}$ \\
$\mathcal C_{13}$ & $\patset{12/RB,21/GY}$,\allowbreak\ $\patset{12/RB,21/YG}$,\allowbreak\ $\patset{12/RG,21/BY}$,\allowbreak\ $\patset{12/RG,21/YB}$,\allowbreak\ $\patset{12/RY,21/BG}$,\allowbreak\ $\patset{12/RY,21/GB}$ \\
$\mathcal C_{13}$ & $\patset{12/BR,21/GY}$,\allowbreak\ $\patset{12/BR,21/YG}$,\allowbreak\ $\patset{12/BB,21/RG}$,\allowbreak\ $\patset{12/BB,21/RY}$,\allowbreak\ $\patset{12/BB,21/GR}$,\allowbreak\ $\patset{12/BB,21/GY}$ \\
$\mathcal C_{13}$ & $\patset{12/BB,21/YR}$,\allowbreak\ $\patset{12/BB,21/YG}$,\allowbreak\ $\patset{12/BG,21/RY}$,\allowbreak\ $\patset{12/BG,21/YR}$,\allowbreak\ $\patset{12/BY,21/RG}$,\allowbreak\ $\patset{12/BY,21/GR}$ \\
$\mathcal C_{13}$ & $\patset{12/GR,21/BY}$,\allowbreak\ $\patset{12/GR,21/YB}$,\allowbreak\ $\patset{12/GB,21/RY}$,\allowbreak\ $\patset{12/GB,21/YR}$,\allowbreak\ $\patset{12/GG,21/RB}$,\allowbreak\ $\patset{12/GG,21/RY}$ \\
$\mathcal C_{13}$ & $\patset{12/GG,21/BR}$,\allowbreak\ $\patset{12/GG,21/BY}$,\allowbreak\ $\patset{12/GG,21/YR}$,\allowbreak\ $\patset{12/GG,21/YB}$,\allowbreak\ $\patset{12/GY,21/RB}$,\allowbreak\ $\patset{12/GY,21/BR}$ \\
$\mathcal C_{13}$ & $\patset{12/YR,21/BG}$,\allowbreak\ $\patset{12/YR,21/GB}$,\allowbreak\ $\patset{12/YB,21/RG}$,\allowbreak\ $\patset{12/YB,21/GR}$,\allowbreak\ $\patset{12/YG,21/RB}$,\allowbreak\ $\patset{12/YG,21/BR}$ \\
$\mathcal C_{13}$ & $\patset{12/YY,21/RB}$,\allowbreak\ $\patset{12/YY,21/RG}$,\allowbreak\ $\patset{12/YY,21/BR}$,\allowbreak\ $\patset{12/YY,21/BG}$,\allowbreak\ $\patset{12/YY,21/GR}$,\allowbreak\ $\patset{12/YY,21/GB}$ \\
$\mathcal C_{13}$ & $\patset{21/RR,21/BG}$,\allowbreak\ $\patset{21/RR,21/BY}$,\allowbreak\ $\patset{21/RR,21/GB}$,\allowbreak\ $\patset{21/RR,21/GY}$,\allowbreak\ $\patset{21/RR,21/YB}$,\allowbreak\ $\patset{21/RR,21/YG}$ \\
$\mathcal C_{13}$ & $\patset{21/RB,21/GG}$,\allowbreak\ $\patset{21/RB,21/YY}$,\allowbreak\ $\patset{21/RG,21/BB}$,\allowbreak\ $\patset{21/RG,21/YY}$,\allowbreak\ $\patset{21/RY,21/BB}$,\allowbreak\ $\patset{21/RY,21/GG}$ \\
$\mathcal C_{13}$ & $\patset{21/BR,21/GG}$,\allowbreak\ $\patset{21/BR,21/YY}$,\allowbreak\ $\patset{21/BB,21/GR}$,\allowbreak\ $\patset{21/BB,21/GY}$,\allowbreak\ $\patset{21/BB,21/YR}$,\allowbreak\ $\patset{21/BB,21/YG}$ \\
$\mathcal C_{13}$ & $\patset{21/BG,21/YY}$,\allowbreak\ $\patset{21/BY,21/GG}$,\allowbreak\ $\patset{21/GR,21/YY}$,\allowbreak\ $\patset{21/GB,21/YY}$,\allowbreak\ $\patset{21/GG,21/YR}$,\allowbreak\ $\patset{21/GG,21/YB}$ \\
\midrule
$\mathcal C_{14}$ & $\patset{12/RB,21/BG}$,\allowbreak\ $\patset{12/RB,21/BY}$,\allowbreak\ $\patset{12/RG,21/GB}$,\allowbreak\ $\patset{12/RG,21/GY}$,\allowbreak\ $\patset{12/RY,21/YB}$,\allowbreak\ $\patset{12/RY,21/YG}$ \\
$\mathcal C_{14}$ & $\patset{12/BR,21/RG}$,\allowbreak\ $\patset{12/BR,21/RY}$,\allowbreak\ $\patset{12/BG,21/GR}$,\allowbreak\ $\patset{12/BG,21/GY}$,\allowbreak\ $\patset{12/BY,21/YR}$,\allowbreak\ $\patset{12/BY,21/YG}$ \\
$\mathcal C_{14}$ & $\patset{12/GR,21/RB}$,\allowbreak\ $\patset{12/GR,21/RY}$,\allowbreak\ $\patset{12/GB,21/BR}$,\allowbreak\ $\patset{12/GB,21/BY}$,\allowbreak\ $\patset{12/GY,21/YR}$,\allowbreak\ $\patset{12/GY,21/YB}$ \\
$\mathcal C_{14}$ & $\patset{12/YR,21/RB}$,\allowbreak\ $\patset{12/YR,21/RG}$,\allowbreak\ $\patset{12/YB,21/BR}$,\allowbreak\ $\patset{12/YB,21/BG}$,\allowbreak\ $\patset{12/YG,21/GR}$,\allowbreak\ $\patset{12/YG,21/GB}$ \\
\midrule
$\mathcal C_{15}$ & $\patset{12/RB,21/RG}$,\allowbreak\ $\patset{12/RB,21/RY}$,\allowbreak\ $\patset{12/RG,21/RB}$,\allowbreak\ $\patset{12/RG,21/RY}$,\allowbreak\ $\patset{12/RY,21/RB}$,\allowbreak\ $\patset{12/RY,21/RG}$ \\
$\mathcal C_{15}$ & $\patset{12/BR,21/BG}$,\allowbreak\ $\patset{12/BR,21/BY}$,\allowbreak\ $\patset{12/BG,21/BR}$,\allowbreak\ $\patset{12/BG,21/BY}$,\allowbreak\ $\patset{12/BY,21/BR}$,\allowbreak\ $\patset{12/BY,21/BG}$ \\
$\mathcal C_{15}$ & $\patset{12/GR,21/GB}$,\allowbreak\ $\patset{12/GR,21/GY}$,\allowbreak\ $\patset{12/GB,21/GR}$,\allowbreak\ $\patset{12/GB,21/GY}$,\allowbreak\ $\patset{12/GY,21/GR}$,\allowbreak\ $\patset{12/GY,21/GB}$ \\
$\mathcal C_{15}$ & $\patset{12/YR,21/YB}$,\allowbreak\ $\patset{12/YR,21/YG}$,\allowbreak\ $\patset{12/YB,21/YR}$,\allowbreak\ $\patset{12/YB,21/YG}$,\allowbreak\ $\patset{12/YG,21/YR}$,\allowbreak\ $\patset{12/YG,21/YB}$ \\
\midrule
$\mathcal C_{16}$ & $\patset{21/RB,21/BG}$,\allowbreak\ $\patset{21/RB,21/BY}$,\allowbreak\ $\patset{21/RB,21/GR}$,\allowbreak\ $\patset{21/RB,21/YR}$,\allowbreak\ $\patset{21/RG,21/BR}$,\allowbreak\ $\patset{21/RG,21/GB}$ \\
$\mathcal C_{16}$ & $\patset{21/RG,21/GY}$,\allowbreak\ $\patset{21/RG,21/YR}$,\allowbreak\ $\patset{21/RY,21/BR}$,\allowbreak\ $\patset{21/RY,21/GR}$,\allowbreak\ $\patset{21/RY,21/YB}$,\allowbreak\ $\patset{21/RY,21/YG}$ \\
$\mathcal C_{16}$ & $\patset{21/BR,21/GB}$,\allowbreak\ $\patset{21/BR,21/YB}$,\allowbreak\ $\patset{21/BG,21/GR}$,\allowbreak\ $\patset{21/BG,21/GY}$,\allowbreak\ $\patset{21/BG,21/YB}$,\allowbreak\ $\patset{21/BY,21/GB}$ \\
$\mathcal C_{16}$ & $\patset{21/BY,21/YR}$,\allowbreak\ $\patset{21/BY,21/YG}$,\allowbreak\ $\patset{21/GR,21/YG}$,\allowbreak\ $\patset{21/GB,21/YG}$,\allowbreak\ $\patset{21/GY,21/YR}$,\allowbreak\ $\patset{21/GY,21/YB}$ \\
\midrule
$\mathcal C_{17}$ & $\patset{21/RB,21/GB}$,\allowbreak\ $\patset{21/RB,21/YB}$,\allowbreak\ $\patset{21/RG,21/BG}$,\allowbreak\ $\patset{21/RG,21/YG}$,\allowbreak\ $\patset{21/RY,21/BY}$,\allowbreak\ $\patset{21/RY,21/GY}$ \\
$\mathcal C_{17}$ & $\patset{21/BR,21/GR}$,\allowbreak\ $\patset{21/BR,21/YR}$,\allowbreak\ $\patset{21/BG,21/YG}$,\allowbreak\ $\patset{21/BY,21/GY}$,\allowbreak\ $\patset{21/GR,21/YR}$,\allowbreak\ $\patset{21/GB,21/YB}$ \\
\midrule
$\mathcal C_{18}$ & $\patset{21/RB,21/RG}$,\allowbreak\ $\patset{21/RB,21/RY}$,\allowbreak\ $\patset{21/RG,21/RY}$,\allowbreak\ $\patset{21/BR,21/BG}$,\allowbreak\ $\patset{21/BR,21/BY}$,\allowbreak\ $\patset{21/BG,21/BY}$ \\
$\mathcal C_{18}$ & $\patset{21/GR,21/GB}$,\allowbreak\ $\patset{21/GR,21/GY}$,\allowbreak\ $\patset{21/GB,21/GY}$,\allowbreak\ $\patset{21/YR,21/YB}$,\allowbreak\ $\patset{21/YR,21/YG}$,\allowbreak\ $\patset{21/YB,21/YG}$ \\
\midrule
$\mathcal C_{19}$ & $\patset{21/RB,21/GY}$,\allowbreak\ $\patset{21/RB,21/YG}$,\allowbreak\ $\patset{21/RG,21/BY}$,\allowbreak\ $\patset{21/RG,21/YB}$,\allowbreak\ $\patset{21/RY,21/BG}$,\allowbreak\ $\patset{21/RY,21/GB}$ \\
$\mathcal C_{19}$ & $\patset{21/BR,21/GY}$,\allowbreak\ $\patset{21/BR,21/YG}$,\allowbreak\ $\patset{21/BG,21/YR}$,\allowbreak\ $\patset{21/BY,21/GR}$,\allowbreak\ $\patset{21/GR,21/YB}$,\allowbreak\ $\patset{21/GB,21/YR}$ \\
\end{longtable}
\endgroup

\end{document}